\documentclass[singlespacing]{elsart}

\usepackage[usenames]{color}
\usepackage{epsfig,graphicx}
\usepackage{psfrag}
\usepackage{natbib}
\usepackage{epsfig}
\usepackage{amsmath,amsfonts,amsxtra,rotating,array,multirow,amsxtra}
\usepackage{amssymb,bm}

\newtheorem{theorem}{Theorem}[section]
\newtheorem{lemma}[theorem]{Lemma}

\newtheorem{corollary}[theorem]{Corollary}

\date{}

\journal {Elsevier}
\begin{document}

\begin{frontmatter}



\title{\color{black}A non-adapted sparse approximation of PDEs with stochastic inputs}


\author[CU]{Alireza Doostan\corauthref{cor}}
 \corauth[cor]{Corresponding author. Tel.:
+(303)-492-7572; Fax: +(303)-492-4990.}
\ead{doostan@colorado.edu}
\author[CIT]{Houman Owhadi}
\ead{owhadi@caltech.edu}
\address[CU]{Aerospace Engineering Sciences Department, University of Colorado, Boulder, CO 80309, USA}
\address[CIT]{Applied \& Computational Mathematics Department, California Institute of Technology, Pasadena, CA 91125, USA}

\begin{abstract}
We propose a method for the approximation of solutions of PDEs with stochastic coefficients based on the direct, i.e., non-adapted, sampling of solutions. This sampling can be done by using any legacy code for the deterministic problem as a black box. The method converges in probability (with probabilistic error bounds) as a consequence of sparsity and a concentration of measure phenomenon on the empirical correlation between samples. We show that the method is well suited for truly high-dimensional problems (with slow decay in the spectrum). 
\end{abstract}

\begin{keyword}
{\color{black}Polynomial chaos; Uncertainty quantification; Stochastic PDE; Compressive sampling; Sparse approximation}
\end{keyword}

\end{frontmatter}
\section{Introduction}
\label{sec:intro}

Realistic analysis and design of complex engineering systems require not only a fine understanding and modeling of the underlying physics and their interactions, but also a significant recognition of intrinsic uncertainties and their influences on the quantities of interest. Uncertainty Quantification (UQ) is an emerging discipline that aims at addressing the latter issue; it aims at meaningful {characterization} of uncertainties in the physical models from the available measurements and efficient propagation of these uncertainties for a quantitative validation of model predictions.

Despite recent growing interests in UQ of complex systems, it remains a grand challenge to efficiently propagate uncertainties through systems characterized by a large number of uncertain sources where the so-called curse-of-dimensionality is yet an unsolved problem. Additionally, development of {\it non-intrusive} uncertainty propagation techniques is of essence as the analysis of complex multi-disciplinary systems often requires the use of sophisticated coupled deterministic solvers in which one cannot readily intrude to set up the necessary propagation infrastructure.

Sampling methods such as the Monte Carlo simulation and its several variants had been utilized for a long time as the primary scheme for uncertainty propagation. However, it is well understood that these methods are generally inefficient for large-scale systems due to their slow rate of convergence. There has been an increasing recent interest in developing alternative numerical methods that are more efficient than the Monte Carlo techniques. Most notably, the {stochastic Galerkin} schemes using polynomial chaos (PC) bases \citep{Ghanem03,Deb01,Xiu02,Babuska04,Wan05} have been successfully applied to a variety of engineering problems and are extremely useful when the number of uncertain parameters is not large. In their original form, the stochastic Galerkin schemes are {\it intrusive}, as one has to modify the deterministic solvers for their implementation. {\it Stochastic collocation} schemes \citep{Tatang95,Mathelin03,Xiu05a,Babuska07a, Nobile08a} belong to a different class of methods that rely upon (isotropic) sparse grid integration/interpolation in the stochastic space of the problem to reduce the curse-of-dimensionality associated with the conventional tensor-product integration/interpolation rules. As their construction is primarily based on the input parameter space, the computational cost of both stochastic Galerkin and collocation techniques increases rapidly for large number of independent input uncertainties.

More recently, efforts have been made to construct {\it solution-adaptive} uncertainty propagation techniques that exploit any structures in the solution to decrease the computational cost. Among them are the multi-scale model reduction of  \citep{Doostan07} and the sparse decomposition of \citep{Todor07a,Bieri09a, Bieri09b,Bieri09c,Blatman10} for the stochastic Galerkin technique, anisotropic and adaptive sparse grids of \citep{Nobile08b,Ma09a} for the stochastic collocation scheme, and low-rank solution approximations of \citep{Nouy07,Nouy08,Doostan09}.

In the present study, we are interested in cases where the quantity of interest is {\it sparse} at the stochastic level, i.e., it can be accurately represented with only {\it few} terms when linearly expanded into a stochastic, e.g., polynomial chaos, basis. Interestingly, sparsity is salient in the analysis of high-dimensional problems where the number of energetic basis functions (those with large coefficients) is small relative to the cardinality of the full basis. For instance, it has been shown in \citep{Todor07a,Bieri09a} that, under some mild conditions, solutions to linear elliptic stochastic PDEs with high-dimensional random coefficients admit sparse representations with respect to the PC basis. Consequently, an approach based on a zero-dimensional algebraic stochastic problem has been proposed in \citep{Bieri09a} to detect the sparsity pattern, which then guides the stochastic Galerkin analysis of the original problem. Moreover, a ``quasi"-best $N$-term approximation for a class of elliptic stochastic PDEs has been proposed in \citep{Bieri09c}.

In this work, using  {\it concentration of measure} inequalities and  {\it compressive sampling} techniques, we derive a method for PC expansion of sparse solutions to stochastic PDEs. The proposed method is

 \begin{itemize}
 \item Non-intrusive: it is based on the direct random sampling of the PDE solutions. This sampling can be done by using any legacy code for the deterministic problem as a black box.
 \item Non-adapted: it does not tailor the sampling process to identify the important dimensions at the stochastic level
 \item Provably convergent: we obtain probabilistic bounds on the approximation error proving the stability and convergence of the method.
 \item {\color{black}Well-suited to problems with high-dimensional random inputs.}
 \end{itemize}

{\it Compressive sampling} is an emerging direction in signal processing that aims at recovering {\it sparse} signals accurately (or even exactly) from a small number of their random projections \citep{Chen98,Chen01a,Candes06a,Donoho06b,Candes06b,Candes06c,Candes07a,Cohen09a,Bruckstein09}. A sparse signal is simply a signal that has only few significant coefficients when linearly expanded into a basis, e.g., $\{\psi_{\bm{\alpha}}\}$.

For sufficiently sparse signals, the number of samples needed for a successful recovery is typically less than what is required by the Shannon-Nyquist sampling principle. Generally speaking, a successful signal reconstruction by compressive sampling is conditioned upon:
\begin{itemize}
\item {\it Sufficient} sparsity of the signal; and
\item {\it Incoherent} random projections of the signal.
\end{itemize}
A square-measurable stochastic function $u(\omega)$, defined on a suitable probability space $(\Omega,\mathcal{F},\mathcal{P})$ can be expanded into a mean-squared convergent series of the chaos polynomial bases, i.e., $u(\omega)\approx\sum_{\bm \alpha}c_{\bm \alpha}\psi_{\bm \alpha}(\omega)$, with some cardinality $P$. The stochastic function $u(\omega)$ is then sparse in PC basis $\{\psi_{\bm \alpha}\}$, if only a small fraction of coefficients $c_{\bm \alpha}$ are significant. In this case, under certain conditions, the sparse PC coefficients $\bm{c}$ may be computed accurately and robustly using only $N\ll P$ random samples of $u(\omega)$ via compressive sampling. Given $N$ random samples of $u(\omega)$, compressive sampling aims at finding the sparsest (or nearly sparsest) coefficients $\bm{c}$ from an optimization problem of the form
\begin{equation}
\label{eqn:P1_delta_intro}
(P_{s,\delta}):\qquad\min_{\bm{c}} \Vert\bm{Wc}\Vert_{s}\quad\mathrm{subject \;to}\quad \Vert\bm{\Psi}\bm{c}-\bm{u}\Vert_2\leq\delta,
\end{equation}
\noindent where $\Vert\bm{Wc}\Vert_{s}$, with $s=\{0,1\}$ and some positive diagonal weight matrix $\bm W$, is a measure of the sparsity of $\bm{c}$ and $\Vert\bm{\Psi}\bm{c}-\bm{u}\Vert_2$ is a measure of the accuracy of the truncated PC expansion in estimating the $u(\omega)$ samples. The $N$-vector $\bm{u}$ contains the independent random samples of $u(\omega)$ and the rows of the $N\times P$ matrix $\bm{\Psi}$ consist of the corresponding samples of the PC basis $\{\psi_{\bm \alpha}\}$.

Throughout the rest of this manuscript, we will elaborate on the formulation of the compressive sampling problem (\ref{eqn:P1_delta_intro}) and the required conditions under which it leads to an accurate and stable approximation of an arbitrary sparse stochastic function as well as sparse solutions to linear elliptic stochastic PDEs. {\color{black}Although we choose to study this particular class of stochastic PDEs, we stress that the proposed algorithms and theoretical developments are far more general and may be readily applied to recover sparse solution of other stochastic systems.}

{\color{black} In Section \ref{sec:p-setup}, we describe the setup of the problem of interest. We then, in Section \ref{sec:stoch_disc}, briefly overview the spectral stochastic discretization of the random functions using PC basis. The main contribution of this work on sparse approximation of stochastic PDEs using compressive sampling is then introduced in Sections \ref{sec:cs}, \ref{sec:stability_spde}, and \ref{sec:sucess_rec}. Sections \ref{sec:delta} and \ref{sec:solvers} discuss some of the implementation details of the present technique. To demonstrate the accuracy and efficiency of the proposed procedures, in Section \ref{sec:examples}, we perform two numerical experiments on a $1$-$D$ (in space) linear elliptic stochastic PDE with high-dimensional random diffusion coefficients.}

{\color{black}
\section{Problem setup}
\label{sec:p-setup}}
Let $(\Omega,\mathcal{F},\mathcal{P})$ be a complete probability space where $\mathcal{P}$ is a probability measure on the   $\sigma-$field  $\mathcal{F}$. We consider the following elliptic stochastic PDE defined on a bounded {\color{black}Lipschitz continuous} domain $\mathcal{D}\subset \mathbb{R}^{D}$, $D=1,2,$ or $3$, with boundary $\partial \mathcal{D}$,
\begin{eqnarray}
\label{eqn:spde}
-\nabla\cdot\left(a(\bm{x},\omega)\nabla u(\bm{x},\omega)\right)&=&f(\bm{x})\quad \bm{x}\in \mathcal{D},\\
u(\bm{x},\omega)&=&0\quad \bm{x}\in \partial\mathcal{D}, \nonumber
\end{eqnarray}
\noindent $\mathcal{P}-a.s.\; \omega\in\Omega$. The diffusion coefficient $a(\bm{x},\omega)$ is a stochastic function defined on $(\Omega,\mathcal{F},\mathcal{P})$ and is the source of uncertainty in (\ref{eqn:spde}). We assume that $a(\bm{x},\omega)$ is specified by a truncated Karhunen-Lo\`eve-``like" expansion
\begin{equation}
\label{eqn:kle_a}
a(\bm{x},\omega)=\bar{a}(\bm{x})+\sum_{i=1}^{d}\sqrt{\lambda_i}\phi_{i}(\bm{x})y_{i}(\omega),
\end{equation}
\noindent where $(\lambda_i,\phi_i)$, $i=1,\cdots,d$,  are the eigenpairs of the covariance function $C_{aa}(\bm{x}_1,\bm{x}_2)\in L_{2}(\mathcal{D}\times\mathcal{D})$ of $a(\bm{x},\omega)$ and $\bar{a}(\bm{x})$ is the mean of $a(\bm{x},\omega)$. We further assume that $a(\bm{x},\omega)$ satisfies the following conditions:

\begin{description}

\item {\it A-I.} For all $\bm{x}\in \mathcal{D}$, there exists constants $a_{\min}$ and $a_{\max}$ such that
\begin{equation}
\label{eqn:a_positivity}
 0<a_{\min}\leq a(\bm{x},\omega)\leq a_{\max}<\infty \quad \mathcal{P}-a.s.\; \omega\in\Omega
\end{equation}
\item {\it A-II.} The covariance function $C_{aa}(\bm{x}_1,\bm{x}_2)$ is {\it piecewise analytic} on $\mathcal{D}\times\mathcal{D}$  {\color{black}\citep{Schwab06a,Bieri09a}}, implying that there exist real constants $c_1$ and $c_2$ such that for $i=1,\cdots,d$,
\begin{equation}
\label{eqn:lambda_bound}
0\leq \lambda_i\leq c_1e^{-c_2i^{\kappa}}
\end{equation}
and
\begin{equation}
\label{eqn:phi_bound}
\forall \bm{\alpha}\in\mathbb{N}^{d}:\qquad \sqrt{\lambda_{i}}\Vert\partial^{\bm{\alpha}}\phi_i\Vert_{L^{\infty}(\mathcal{D})}\leq c_1e^{-c_2i^{\kappa}},
\end{equation}
\noindent where $\kappa:=1/D$ and  $\bm{\alpha}\in\mathbb{N}^{d}$ is a fixed multi-index. Notice that the decay rates in Eqs. (\ref{eqn:lambda_bound}) and (\ref{eqn:phi_bound}) will be {\it algebraic} if $C_{aa}(\bm{x}_1,\bm{x}_2)$ has $C^{s}(\mathcal{D}\times\mathcal{D})$ regularity for some $s>0$ {\color{black}\citep{Schwab06a}}. \\

\item {\it A-III.} The random variables $\{y_{k}(\omega)\}_{k=1}^{d}$ are independent and uniformly distributed on $\Gamma_k:=[-1,1]$, $k=1,\cdots,d$, with probability distribution function $\rho_k(y_k)=1/2$ defined over $\Gamma_k$. The joint probability distribution function of the random vector {\color{black}$\bm{y}:=(y_{1},\cdots,y_d)$} is then given by $\rho(\bm{y}):=\prod_{k=1}^{d}\rho_{k}(y_k)$.\\\

\end{description}

{\bf Remark:} The algorithm proposed in the paper requires the existence of a sparse solution. The only role of
assumption {\it A-II} is  to guarantee the existence of such a sparse solution. It is not necessary for the application and the validity of the proposed algorithm. In particular, if the coefficient $a$ is only essentially bounded, the proposed algorithm will be accurate as long as a sparse approximation exists.

Given the {\color{black}finite-dimensional uncertainty representation in (\ref{eqn:kle_a})}, the solution $u(\bm{x},\omega)$ of (\ref{eqn:spde}) also admits a finite-dimensional representation, i.e.,
\begin{equation}
\label{eqn:sol}
u(\bm{x},\bm{y}):=u(\bm{x},y_1(\omega),\cdots,y_d(\omega)):\mathcal{D}\times\Gamma\rightarrow \mathbb{R},
\end{equation}
\noindent where $\Gamma:=\prod_{k=1}^{d}\Gamma_{k}$.

{\color{black}In what follows, we first briefly outline the Legendre spectral stochastic discretization of $u(\bm{x},\bm{y})$ and consequently introduce our approach based on compressive sampling to obtain such a discretization.}
{\color{black}
\section{Numerical approach}
\label{sec:numerical_approach}}
\subsection{Spectral stochastic discretization}
\label{sec:stoch_disc}

{\color{black}In the context of the spectral stochastic methods {\color{black}\citep{Ghanem03,Deb01,Xiu02,Babuska04}}, the solution $u(\bm{x},\bm{y})$ of (\ref{eqn:spde}) is represented by an infinite series of the form
\begin{equation}
\label{eqn:chaos_exact}
u(\bm{x},\bm{y})=\sum_{\bm{\alpha}\in\mathbb{N}_{0}^{d}}c_{\bm{\alpha}}(\bm{x})\psi_{\bm{\alpha}}(\bm{y}),
\end{equation}
where $\mathbb{N}_{0}^{d}:=\left\{(\alpha_1,\cdots,\alpha_d): \alpha_j\in \mathbb{N}\cup\{0\} \right\}$ is the set of multi-indices of size $d$ defined on non-negative integers.} The basis functions $\{\psi_{\bm{\alpha}}(\bm{y})\}$ are multi-dimensional Legendre polynomials, referred to as the Legendre polynomial chaos (PC), and are orthogonal with respect to the joint probability measure $\rho(\bm{y})$ of the random vector $\bm{y}$. Each basis function $\psi_{\bm{\alpha}}(\bm{y})$ is a tensor product of univariate Legendre polynomials $\psi_{\alpha_i}(y_i)$ of degree $\alpha_i\in\mathbb{N}_0^1$, i.e.,
\begin{equation}
\label{eqn:multid-univ}
\psi_{\bm{\alpha}}(\bm{y})=\psi_{\alpha_1}(y_1)\psi_{\alpha_2}(y_2)\cdots\psi_{\alpha_d}(y_d)\qquad {\bm{\alpha}}\in\mathbb{N}_{0}^{d}.
\end{equation}
We here assume that the univariate Legendre polynomials  $\psi_{\alpha_i}(y_i)$ are also normalized {\color{black}such that}
\begin{equation}
\label{eqn:univ_legendre}
\int_{\Gamma_i}\psi_{\alpha_i}^{2}(y_i)\rho_{i}(y_i)dy_i=1,\quad i=1,\cdots,d.
\end{equation}
The {\color{black}exact} generalized Fourier coefficients $c_{\bm{\alpha}}(\bm{x})$ in (\ref{eqn:chaos_exact}), referred to as the PC coefficients, are computed by the projection of $u(\bm{x},\bm{y})$ onto each basis function $\psi_{\bm{\alpha}}(\bm{y})$,
\begin{equation}
\label{eqn:coefficients}
c_{\bm{\alpha}}(\bm{x})=\mathbb{E}\left[u(\bm{x},\bm{y})\psi_{\bm{\alpha}}(\bm{y})\right] =\int_{\Gamma}u(\bm{x},\bm{y})\psi_{\bm{\alpha}}(\bm{y})\rho(\bm{y})d\bm{y}.
\end{equation}
{\color{black}Here, $\mathbb{E}$ denotes the expectation operator.} In practice, the expansion (\ref{eqn:chaos_exact}) is finite; that is, only a finite number of spectral modes is needed to approximate $u(\bm{x},\bm{y})$ within a given target accuracy. Traditionally, a finite-order truncation of the basis $\{\psi_{\bm{\alpha}}(\bm{y})\}$ is adopted for the approximation, i.e.,
\begin{equation}
\label{eqn:chaos_truncated}
u_{p}(\bm{x},\bm{y}):=\sum_{\bm{\alpha}\in\Lambda_{p,d}}c_{\bm{\alpha}}(\bm{x})\psi_{\bm{\alpha}}(\bm{y}),
\end{equation}
\noindent where the set of multi-indices $\Lambda_{p,d}$ is
\begin{equation}
\label{eqn:multi-index}
\Lambda_{p,d}:=\left\{\bm{\alpha}\in\mathbb{N}_{0}^{d}: \Vert\bm{\alpha}\Vert_{1}\leq p,\;\Vert \bm{\alpha}\Vert_{0}\leq d \right\}
\end{equation}
\noindent and has the cardinality
\begin{equation}
\label{eqn:card_basis_full}
P:=\vert\Lambda_{p,d}\vert=\frac{(p+d)!}{p!d!}.
\end{equation}
Here, $\Vert\bm{\alpha}\Vert_{1}=\sum_{i=1}^{d}\bm{\alpha}_i$ and $\Vert\bm{\alpha}\Vert_{0}=\#\{i:\alpha_i>0\}$ are the total order (degree) and dimensionality of the basis function $\psi_{\bm{\alpha}}(\bm{y})$, respectively. The approximation is then refined by increasing $p$ to achieve a given target accuracy. Under assumptions {\it A-I}, {\it A-II}, and {\it A-III} stated in Section \ref{sec:p-setup}, the solution $u(\bm{x},\bm{y})$ is analytic with respect to the random variables $\{y_i\}_{i=1}^{d}$ (see \citep{Babuska07a}), and as $p$ increases, the approximation (\ref{eqn:chaos_truncated}) converges exponentially fast in the mean-squares sense \citep{Babuska04,Babuska07a,Bieri09a}.

{\bf Definition (Sparsity)} {\it The solution $u(\bm{x},\bm{y})$ is said to be (nearly) sparse if only a small fraction of coefficients $c_{\bm\alpha}(\bm{x})$ in (\ref{eqn:chaos_truncated}) are dominant and contribute to the solution statistics.}

As will be described in Section \ref{sec:cs}, a sparse solution $u(\bm{x},\bm{y})$ may be accurately recovered using $N\ll P$ random samples $\{u(\bm{x},\bm{y}_i)\}_{i=1}^{N}$ using compressive sampling techniques. This has to be compared, for instance, with the least-squares regression-type techniques, \citep{Hosder06}, that normally require $N\gg P$ samples for an accurate recovery.
\subsection{Sparse recovery using compressive sampling}
\label{sec:cs}
{\it Compressive sampling} is an emerging theory in the field of signal and image processing \citep{Chen98,Chen01a,Candes06a,Donoho06b,Candes06b,Candes06c,Candes07a,Cohen09a,Bruckstein09}. It hinges around the idea that a set of incomplete random observations of a sparse signal can be used to accurately, or even exactly, recover the signal (provided that the basis in which the signal is sparse is known). In particular, the number of such observations may be much smaller than the cardinality of the signal. In the context of problem (\ref{eqn:spde}), compressive sampling may be interpreted as follows. The solution $u(\bm{x},\bm{y})$ that is sparse, in the sense of Lemma \ref{lem:sparsity} defined in Section \ref{sec:sparse_solution}, can be accurately recovered using $N\ll P$ random samples $\{u(\bm{x},\bm{y}_i)\}_{i=1}^{N}$, where $P$ is the cardinality of the Legendre PC basis $\{\bm{\psi}_{\bm{\alpha}}(\bm{y})\}$. We next elaborate on the above statement and address how such a sparse reconstruction is achieved and under what conditions it is successful.

Let  $\{u(\bm{y}_{i})\}_{i=1}^{N}$ be {\it i.i.d.} random samples of $u(\bm{x},\bm{y})$ for a fixed point $\bm{x}$ in $\mathcal{D}$. For the time being, let us assume that the $p$th-order PC basis $\{\bm{\psi}_{\bm{\alpha}}(\bm{y})\}$ is a complete basis to expand $u(\bm{y})$; we will relax this assumption as we proceed. Given pairs of $\{\bm{y}_{i}\}_{i=1}^{N}$ and $\{u(\bm{y}_{i})\}_{i=1}^{N}$, we write
\begin{equation}
\label{eqn:pce_of_samples}
u(\bm{y}_{i}) = \sum_{\bm{\alpha}\in\Lambda_{p,d}}c_{\bm{\alpha}}\psi_{\bm{\alpha}}(\bm{y}_i),\quad i=1,\cdots,N,
\end{equation}
\noindent or equivalently,
\begin{equation}
\label{eqn:underdetermind_system}
\bm{\Psi}\bm{c}=\bm{u}.
\end{equation}
Here $\bm{c}\in\mathbb{R}^{P}$ is the vector of PC coefficients $c_{\bm{\alpha}}$ to be determined, $\bm{u}\in\mathbb{R}^{N}$ is the vector of samples of $u(\bm{y})$, and each column of the {\it measurement} matrix $\bm{\Psi}\in \mathbb{R}^{N\times P}$ contains samples of the $j$th element of the PC basis, i.e.,
\begin{equation}
\label{eqn:Psi_matrix}
\bm{\Psi}[i,j] = \psi_{\bm{\alpha}_{j}}(\bm{y}_{i}),\quad i=1,\cdots,N,\quad j=1,\cdots,P.
\end{equation}
We are interested in the case that the number $N$ of solution samples is much smaller than the unknown PC coefficients $P$, i.e., $N\ll P$. Without any additional constraints on $\bm{c}$, the underdetermined linear system (\ref{eqn:underdetermind_system}) is ill-posed and, in general, has infinitely many solutions. When $\bm{c}$ is sparse; that is, only a small fraction of the coefficients $c_{\bm{\alpha}}$ are significant, the problem \eqref{eqn:underdetermind_system} may be regularized to ensure a well-posed solution. Such a regularization may be achieved by seeking a solution $\bm{c}$ with the minimum number of non-zeros. This can be formulated in the optimization problem
\begin{equation}
\label{eqn:P0}
(P_0):\qquad\min_{\bm{c}} \Vert\bm{c}\Vert_{0}\quad\mathrm{subject \;to}\quad \bm{\Psi}\bm{c}=\bm{u},
\end{equation}
\noindent where the semi-norm $\Vert \bm{c}\Vert_0:=\#\{\bm{\alpha}:c_{{\bm{\alpha}}}\neq 0\}$ is the number of non-zero components of $\bm{c}$. In general, the global minimum solution of $(P_0)$ is not unique and is NP-hard to compute: the cost of a global search is exponential in $P$. Further developments in compressive sampling resulted in a convex relaxation of problem $(P_0)$ by minimization of the $\ell_1$-norm of the solution $\bm{c}$ instead, i.e.,
\begin{equation}
\label{eqn:P1}
(P_{1}):\qquad\min_{\bm{c}} \Vert\bm{Wc}\Vert_{1}\quad\mathrm{subject \;to}\quad \bm{\Psi}\bm{c}=\bm{u},
\end{equation}
\noindent where $\bm{W}$ is a diagonal matrix whose $[j,j]$ entry is the $\ell_2$-norm of the $j$th column of $\bm{\Psi}$ and $\Vert\cdot\Vert_{1}$ denotes the $\ell_1$-norm. Notice that the $\ell_1$-norm is the closest convex function to the $\ell_0$-norm that compels the small coefficients $c_{\bm\alpha}$ to be zero, thus promoting the sparsity in the solution. The purpose of weighting the $\ell_1$ cost function with $\bm{W}$ is to prevent the optimization from biasing toward the non-zero entries in $\bm{c}$ whose corresponding columns in $\bm{\Psi}$ have large norms. The problem $(P_1)$ is entitled Basis Pursuit (BP) \citep{Chen98} and its solution can be obtained by linear programming. Since the $\ell_1$-norm functional $\Vert\bm{c}\Vert_{1}$ is convex, the optimization problem $(P_{1})$ admits a unique solution that coincides with the unique solution to problem $(P_0)$ for sufficiently sparse $\bm{c}$ with some constraints on the measurement matrix $\bm\Psi$; e.g., see \citep{Bruckstein09}. 

In general, the $p$th-order PC basis is not complete for the exact representation of $u(\bm{y})$; therefore, we have to account for the truncation error. This can be accommodated in $(P_0)$ and $(P_1)$ by allowing a non-zero residual in the constraint $\bm{\Psi c}=\bm{u}$. Therefore, as in Sections 3.2.1 and 3.2.3 of \citep{Bruckstein09}, the proposed algorithms in this paper are error-tolerant versions of $(P_0)$ and $(P_1)$, with error tolerance $\delta$, i.e.,
\begin{equation}
\label{eqn:P0_delta}
(P_{0,\delta}):\qquad\min_{\bm{c}} \Vert\bm{c}\Vert_{0}\quad\mathrm{subject \;to}\quad \Vert\bm{\Psi}\bm{c}-\bm{u}\Vert_2\leq\delta
\end{equation}
and
\begin{equation}
\label{eqn:P1_delta}
(P_{1,\delta}):\qquad\min_{\bm{c}} \Vert\bm{Wc}\Vert_{1}\quad\mathrm{subject \;to}\quad \Vert\bm{\Psi}\bm{c}-\bm{u}\Vert_2\leq\delta,
\end{equation}
respectively. The latter problem is named Basis Pursuit Denoising (BPDN) in \citep{Chen98} and may be solved using techniques from quadratic programming. We leave the discussion on the available algorithms for solving problems $(P_{1,\delta})$ and $(P_{0,\delta})$ to Section \ref{sec:solvers}. Instead, we henceforth delineate on sufficient conditions under which the BPDN problem $(P_{1,\delta})$ leads to a successful Legendre PC expansion of a general essentially bounded sparse stochastic function $u(\bm y)$ and, subsequently, the sparse solution $u(\bm x,\bm y)$ to the problem \eqref{eqn:spde}. Our results are extensions of those in \citep{Donoho06a,Bruckstein09}, adapted to the case where the measurement matrix $\bm{\Psi}$ consists of random evaluations of the Legendre PC basis $\{\psi_{\bm{\alpha}}\}$. With slight differences that will be remarked accordingly, similar results hold for the case of the $(P_{0,\delta})$ problem.

\begin{theorem}[General stability of $(P_{1,\delta})$]
\label{the:main0}
Let $u(\bm y)$ be an essentially bounded function of i.i.d. random variables $\bm y:=(y_1,\cdots,y_d)$ uniformly distributed on $\Gamma:=[-1,1]^d$. Define
\begin{equation}
\label{eqn:maximum_sparsity}
S_{\max}:=\frac{N}{64 P^{4c_{p,d}}(\ln P)},
\end{equation}
 with
 \begin{equation}
\label{eqn:cpd}
c_{p,d}:=\frac{\ln 3}{2}\frac{p}{\ln\left(\frac{(p+d)!}{p!d!}\right)}.
\end{equation}
 Let $u_p^{1,\delta}(\bm{y}):=\sum_{\bm \alpha\in\Lambda_{p,d}}c_{\bm\alpha}^{1,\delta}\psi_{\bm\alpha}(\bm y)$ be the Legendre PC approximation of $u(\bm{y})$ with coefficients $\bm{c}^{1,\delta}$ computed from the $\ell_1$-minimization problem $(P_{1,\delta})$ with $\delta\geq 0$.
 If there exists a Legendre PC expansion $u_p^{0}(\bm y):=\sum_{\bm\alpha\in\Lambda_{p,d}^{\epsilon}}c_{\bm\alpha}^{0}\psi_{\bm\alpha}(\bm y)$, for some index set $\Lambda_{p,d}^{\epsilon}\subseteq\Lambda_{p,d}$ such that $\left\Vert u-u_{p}^{0}\right\Vert_{L^{\infty}(\Gamma)}\leq \epsilon$ and
  \begin{equation}
\label{eqn:sparsity_condition}
S < S_{\max},
\end{equation}
with $S:=|\Lambda_{p,d}^{\epsilon}|$, then with probability
 \begin{equation}
\label{eqn:prob_success_general}
Prob_{1}\geq 1- 4P^{2-2S_{\max}}-P^{-8S_{\max}}-P^{-8S_{\max} P^{4c_{p,d}}},
\end{equation}
(on the $N$ samples $\{u(\bm{y}_i)\}_{i=1}^{N}$) and for some constants $c_1$ and $c_2$, the solution $u_p^{1,\delta}$ must obey
\begin{equation}
\label{eqn:stability_general}
\left\Vert u-u_{p}^{1,\delta}\right\Vert_{L^{2}(\Gamma)}\leq c_1\epsilon+c_2\frac{\delta}{\sqrt{N}}.
\end{equation}
\end{theorem}

In simple words, Theorem \ref{the:main0} states that if an essentially bounded stochastic function $u(\bm y)$ admits a sparse Legendre PC expansion, then the $\ell_1$-minimization problem $(P_{1,\delta})$ can accurately recover it from a sufficiently large number of random solution samples. The recovery is stable under the truncation error $\Vert\bm{\Psi}\bm{c}-\bm{u}\Vert_2$ and is within a distance of the exact solution that is proportional to the error tolerance $\delta$. It is worth highlighting that no prior knowledge of the sparsity pattern of the PC coefficients $\bm c$ is needed for an accurate recovery.

{\bf Remark:} Based on the conditions \eqref{eqn:maximum_sparsity} and \eqref{eqn:sparsity_condition}, the number $N$ of random samples has to grow like $P^{4c_{p,d}}\ln P$ and also proportional to the number of dominant coefficients $S=\vert\Lambda_{p,d}^{\epsilon}\vert$. Given any order $p$ of the PC expansion, for sufficiently high-dimensional problems, the constant $c_{p,d}<1/4$ (see Lemma \ref{lem:diameter} and Fig. \ref{fig:c_pd}), thus justifying $N\ll P$. In fact, the conditions \eqref{eqn:maximum_sparsity} and \eqref{eqn:sparsity_condition} are too pessimistic; in practice, the number of random samples required for an accurate recovery is much smaller than the theoretical value in \eqref{eqn:maximum_sparsity}. We will elaborate on this statement in Section \ref{sec:sucess_rec}.

{\bf Remark:} Although the BPDN reconstruction is achieved by minimizing the $\ell_1$-norm of the solution, based on \eqref{eqn:stability_general}, the approximation also converges to the exact solution in the mean-squares sense.

{\bf Remark:} A similar theorem holds for the case of the sparse approximation using the $\ell_0$-minimization problem $(P_{0,\delta})$ in \eqref{eqn:P0_delta}. In this case, the condition \eqref{eqn:maximum_sparsity} has to be replaced with $S_{\max}:=\frac{N}{16 P^{4c_{p,d}}(\ln P)}$ which is, in theory, milder than that of the $(P_{1,\delta})$ problem. The error estimate \eqref{eqn:stability_general} also holds with a larger probability, but with different constants $c_1$ and $c_2$.

\subsection{Sparsity of the solution $u(\bm{x},\bm{y})$}
\label{sec:sparse_solution}

Notice that the accurate recovery of $u(\bm y)$ is conditioned upon the existence of a sparse  PC expansion $u_{p}^{0}$ (see Theorem \ref{the:main0}). In fact, this assumption may not hold for an arbitrary stochastic function $u(\bm y)$, as all the elements of the basis set $\{\psi_{\bm\alpha}(\bm y)\}$ may be important. In this case, our sparse approximation still converges to the actual solution but, perhaps, {\color{black}not using as few as $N\ll P$} random solution samples. 

We will now summarize the results of \citep{Todor07a,Bieri09a} on the sparsity of the Legendre PC expansion of the solution $u(\bm{x},\bm{y})$ to the problem (\ref{eqn:spde}). Alternative to the $p$th-order truncated PC expansion of (\ref{eqn:chaos_truncated}), one may ideally seek a proper index set $\Lambda_{p,d}^{\epsilon}\subseteq \Lambda_{p,d}$, with sufficiently large $p$, such that for a given accuracy $\epsilon$
\begin{equation}
\label{eqn:sparse_set}
\Lambda_{p,d}^{\epsilon}: = \arg\min\left\{\vert \tilde{\Lambda}_{p,d}\vert: \tilde{\Lambda}_{p,d}\subseteq\Lambda_{p,d}, \Vert u-\tilde{u}_{p}\Vert_{H_{0}^{1}\left(\mathcal{D},L^{\infty}(\Gamma)\right)}\leq \epsilon\right\},
\end{equation}
\noindent in which $\tilde{u}_{p}(\bm{x},\bm{y}):=\sum_{\bm{\alpha}\in\tilde{\Lambda}_{p,d}}c_{\bm{\alpha}}(\bm{x})\psi_{\bm{\alpha}}(\bm{y})$. This will lead to the so-called {\it sparse approximation} of $u(\bm{x},\bm{y})$ if
\begin{equation}
\label{eqn:def-sparsity}
\vert\Lambda_{p,d}^{\epsilon}\vert\ll\vert\Lambda_{p,d}\vert=P,
\end{equation}
where $\Lambda_{p,d}$ is defined in \eqref{eqn:multi-index}. Such a reduction in the number of basis functions in (\ref{eqn:def-sparsity}) is possible as, given the accuracy $\epsilon$, the {\it effective} dimensionality $\nu$ of $u(\bm{x},\bm{y})$ in $\Gamma$ is potentially smaller than the apparent dimensionality $d$. More precisely, under assumptions \emph{A-I}, \emph{A-II}, and \emph{A-III} stated in Section \ref{sec:p-setup}, the analyses of \citep{Todor07a,Bieri09a} imply that the discretization of $u(\bm{x},\bm{y})$ using a sparse index set $\Lambda_{p,\nu}$,
\begin{equation}
\label{eqn:multi-index-sparse}
\Lambda_{p,\nu}:=\left\{\bm\alpha\in\mathbb{N}_{0}^{d}: \Vert\bm\alpha\Vert_{1}\leq p,\;\Vert \bm\alpha\Vert_{0}\leq \nu\leq d \right\}
\end{equation}
preserves the exponential decay of the approximation error in the $H_{0}^{1}\left(\mathcal{D},L^{\infty}(\Gamma)\right)$ sense. For the sake of completeness, we cite this from \citep{Bieri09a} in the following lemma.

 \begin{lemma}[Proposition 3.10 of  \citep{Bieri09a}]
 \label{lem:sparsity}
 Given assumptions A-I, A-II, and A-III in Section \ref{sec:p-setup}, there exist constants $c_1,c_2,c_3,c_4 > 0$, depending only on $a(\bm{x},\omega)$ and $f(\bm{x})$ but independent of $d,p,\nu$, such that
\begin{equation}
\label{eqn:err_estimate}
\Vert u-u_{p,\nu}\Vert_{H_{0}^{1}\left(\mathcal{D},L^{\infty}(\Gamma)\right)}\leq c_1\left(e^{-c_2\nu^{1+\kappa}}+e^{c_3\nu(\ln d+\ln p)-c_4 p}\right),
\end{equation}
\noindent {\it for any $d,p,\nu\in\mathbb{N}$ with $\nu\leq d$ and $\kappa=1/D$.}
\end{lemma}
In particular, for $d\geq c_d\vert\ln\epsilon\vert^{1/\kappa}$, choosing
\begin{equation}
\label{eqn:sparse_parameters}
p_\epsilon=\lceil c_pd^{\kappa}\rceil\leq p \quad \mathrm{and}\quad \nu_\epsilon=\lceil c_\nu d^{\kappa/(\kappa+1)}\rceil\leq d,
\end{equation}
\noindent leads to
\begin{equation}
\label{eqn:err_estimate_sparse}
\Vert u-u_{p_\epsilon,\nu_\epsilon}\Vert_{H_{0}^{1}\left(\mathcal{D},L^{\infty}(\Gamma)\right)}\leq\epsilon,
\end{equation}
\noindent where $u_{p_\epsilon,\nu_\epsilon}$ is now defined on a sparse index set
\begin{equation}
\label{eqn:sparse-set-epsilon}
\Lambda_{p_\epsilon,\nu_{\epsilon}}:=\left\{\bm{\alpha}\in\mathbb{N}_{0}^{d}: \Vert\bm\alpha\Vert_{1}\leq p_\epsilon,\;\Vert \bm\alpha\Vert_{0}\leq \nu_\epsilon\right\}
\end{equation}
with cardinality
\begin{equation}
\label{eqn:sparse-set-reduced-dim}
\vert \Lambda_{p_\epsilon,\nu_\epsilon}\vert\lesssim\epsilon^{-1/\rho},
\end{equation}
\noindent for some arbitrary large $\rho>0$ and constants $c_d,c_p$, and $c_\nu$ independent of $d,p_\epsilon$, and $\nu_\epsilon$ \citep{Bieri09a}.

In practice, the sparse set $\Lambda_{p,d}^{\epsilon}$ in (\ref{eqn:sparse_set}) (or equivalently $\Lambda_{p_\epsilon,\nu_\epsilon}$ in (\ref{eqn:sparse-set-epsilon})) is not known {\it a priori}. In \citep{Bieri09a}, an approach based on an algebraic purely-stochastic problem is proposed to adaptively identify $\Lambda_{p,d}^{\epsilon}$. Having done this, the coefficients of the the spectral modes are computed via the (intrusive) stochastic Galerkin scheme \citep{Ghanem03,Xiu02}.
Alternatively, in this work, we apply our sparse approximation using $(P_{1,\delta})$ and $(P_{0,\delta})$ to compute $u(\bm x,\bm y)$. The implementation of $(P_{1,\delta})$ and $(P_{0,\delta})$ is non-intrusive; only random samples of the solution are needed. Moreover, we do not adapt the sampling process to identify the important dimensions at the stochastic level; therefore, our constructions are non-adapted.

Throughout the rest of the present paper, we focus our attention on the case of the stochastic PDE \eqref{eqn:spde} whose solution is provably sparse. The statement of Theorem \ref{the:main0} can be specialized to the approximation of the sparse solution to the stochastic PDE \eqref{eqn:spde} using $(P_{1,\delta})$ (or $(P_{0,\delta})$) as follows.

\subsection{Stability of $(P_{1,\delta})$ for stochastic PDE \eqref{eqn:spde}}
\label{sec:stability_spde}

Combining Lemma \ref{lem:sparsity} with Theorem \ref{the:main0} leads to the following theorem.

\begin{theorem}[Stability of $(P_{1,\delta})$ for stochastic PDE \eqref{eqn:spde}]
\label{the:main}
{\it Let $u_{p}^{1,\delta}(\bm{x},\bm y):=\sum_{\bm\alpha\in\Lambda_{p,d}}c_{\bm\alpha}^{1,\delta}(\bm x)\psi_{\bm\alpha}(\bm y)$ be the $p$th-order Legendre PC approximation of $u(\bm{x},\bm{y})$ in \eqref{eqn:spde} where the coefficients $c_{\bm\alpha}^{1,\delta}(\bm x)$ are obtained from the $\ell_1$-minimization problem $(P_{1,\delta})$ with $N$ independent samples of  $u(\bm{x},\bm{y})$ and arbitrary $\delta\geq 0$.
Write $\kappa:=1/D$. Let $c_{p,d}$ be defined by  \eqref{eqn:cpd} and let $S_{\max}$ be defined by \eqref{eqn:maximum_sparsity}. Let $\rho>0$  be arbitrary.}

{\it There exists constants $c_1,c_2,c_3,c_4,c_5$ independent from $p,d,\kappa,N$ such that 
if  $ \lceil c_2 d^{\kappa}\rceil\leq p$ and $\lceil c_3 d^{\kappa/(\kappa+1)}\rceil\leq d$,}

{\it then with probability at least}
 \begin{equation}
\label{eqn:coherence_boundse2}
Prob_{1}\geq 1- 4P^{2-2S_{\max}}-P^{-8S_{\max}}-P^{-8P^{4c_{p,d}}S_{\max}},
\end{equation}
{\it the solution $u_{p}^{1,\delta}$ must obey}
\begin{equation}
\label{eqn:stability_convergence}
\left\Vert u-u_{p}^{1,\delta}\right\Vert_{L^2(\mathcal{D}, L^2(\Gamma))} \leq c_{4} \epsilon+c_5 \frac{\delta}{\sqrt{N}}.
\end{equation}
with
\begin{equation}
\epsilon:=\max\left(\frac{1}{S_{\max}^\rho},\exp\left(-\left(\frac{d}{c_1}\right)^\kappa\right)\right)
\end{equation}
\end{theorem}

\subsection{Proofs and further ingredients of successful sparse approximations via $(P_{1,\delta})$ and $(P_{0,\delta})$}
\label{sec:sucess_rec}

The ability of problems $(P_{1,\delta})$ and $(P_{0,\delta})$ in accurately approximating the sparse PC coefficients $\bm{c}$ in \eqref{eqn:chaos_truncated}, hence the solution $u(\bm x,\bm{y})$, depends on two main factors: $i)$ the sparsity of the PC coefficients $\bm{c}$ and $ii)$ the {\it mutual coherence} of the measurement matrix $\bm{\Psi}$. In fact, the number $N$ of random solution samples required for a successful sparse approximation is dictated by these two factors. While sparsity is a characteristic of the solution of interest $u(\bm x,\bm y)$, the mutual coherence of the measurement matrix $\bm \Psi$ is universal as it only depends on the choice of PC basis $\{\psi_{\bm\alpha}(\bm y)\}$ and the sampling process from which $\bm\Psi$ is assembled.

In Section \ref{sec:sparse_solution}, based on the analysis of \citep{Bieri09a}, we rationalized the sparsity of $u(\bm x,\bm y)$ with respect to the Legendre PC basis. We now give the definition of the mutual coherence of $\bm{\Psi}$ and discuss its role in our sparse approximation using $(P_{1,\delta})$ and $(P_{0,\delta})$.

\subsubsection{Mutual coherence of $\bm \Psi$}
\label{sec:def_mu}
\textbf{Definition (Mutual Coherence \citep{Donoho06a})} {\it The mutual coherence $\mu(\bm{\Psi})$ of a matrix $\bm{\Psi}\in\mathbb{R}^{N\times P}$ is the maximum of absolute normalized inner-products of its columns. Let $\bm{\psi}_{j}$ and $\bm{\psi}_{k}$ be two distinct columns of $\bm{\Psi}$. Then,}
\begin{equation}
\label{eqn:muA}
\mu(\bm{\Psi}):=\max_{1\leq j,k\leq P,\;j\neq k}\frac{\vert \bm{\psi}_j^{T}\bm{\psi}_k\vert}{\Vert \bm{\psi}_j\Vert_2\Vert \bm{\psi}_k\Vert_2}.
\end{equation}
In plain words, the mutual coherence is a measure of how close to orthogonal a matrix is. Clearly, for any general matrix $\bm\Psi$,
\begin{equation}
\label{eqn:coherence_bound2}
0\leq\mu(\bm{\Psi})\leq 1,
\end{equation}
\noindent where the lower bound is achieved, for instance, by unitary matrices. However, for the case of $N<P$, the mutual coherence $\mu(\bm{\Psi})$ is strictly positive. It is well understood that measurement matrices with smaller mutual coherence have a better ability to recover a sparse solution using compressive sampling techniques, e.g., see Lemma \ref{lem:stability_discrete}. Therefore, we shall proceed to examine the mutual coherence of the random measurement matrix $\bm{\Psi}$ in (\ref{eqn:underdetermind_system}). We first observe that, by the orthogonality of the Legendre PC basis, the mutual coherence $\mu(\bm{\Psi})$ converges to zero almost surely for asymptotically large random sample sizes $N$. However, it is essential for our purpose to $i)$ investigate if a desirably small $\mu(\bm{\Psi})$ can be achieved by a sample size $N\ll P$ and $ii)$ quantify how large $\mu(\bm{\Psi})$ can get for a {\it finite} $N$. These are addressed in the following theorem.

\begin{theorem}[Bound on the mutual coherence $\mu({\bm\Psi})$]
\label{the:mu}
{\it Let $\bm{\Psi}\in\mathbb{R}^{N\times P}$, {\color{black}as defined in (\ref{eqn:underdetermind_system})}, be the measurement matrix corresponding to $N$ independent random samples of the Legendre polynomial chaos basis of order $p$ in $d$ {\it i.i.d.} uniform random variables $\bm{y}$. There exists a positive constant $c_{p,d}:=\frac{\ln 3}{2}\frac{p}{\ln\left(\frac{(p+d)!}{p!d!}\right)}$ depending on $p$ and $d$, such that if}
 \begin{equation}
\label{eqn:define_r}
0\leq r=2\sqrt{\zeta P^{4c_{p,d}}(\ln P)/N}\leq 1/2,
\end{equation}
 {\it for some $\zeta>1$, then}
\begin{equation}
\label{eqn:coherence_bound}
Prob\left[ \mu(\bm{\Psi})\geq \frac{r}{1-r}\right]\leq4P^{2-2\zeta}.
\end{equation}
\end{theorem}

\begin{figure}[htb]
    \centering
      \includegraphics[width=5.0in]{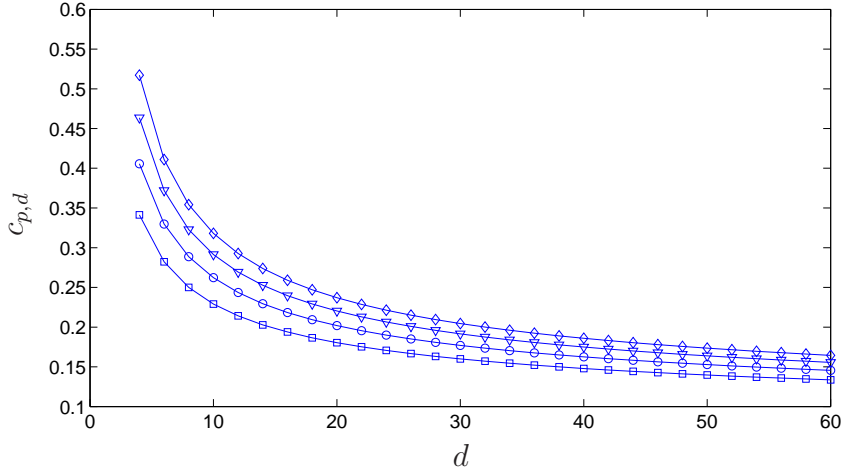}
      \put(-177,-2){$d$}
      \put(-345,85){\begin {sideways} $c_{p,d}$ \end{sideways}}
    \caption{Decay of $c_{p,d}$ as a function of $d$. $p=1$ ({\scriptsize$\square$});   $p=2$ ({$\circ$}); $p=3$ ($\triangledown$);  $p=4$($\diamond$).}
    \label{fig:c_pd}
\end{figure}

Figure \ref{fig:c_pd} illustrates the decay of $c_{p,d}$, for several values of $p$, as a function of $d$. Based on Theorem \ref{the:mu}, for cases where the number $d$ of random variables $\bm{y}$ is large enough such that $c_{p,d}<1/4$, it is sufficient to have $N\sim \mathcal{O}(16P^{4c_{p,d}}\ln P)\ll P$ to keep $\mu(\bm{\Psi})$ bounded from above with a large probability. Notice that such a requirement on $c_{p,d}$ is particularly suited to high-dimensional problems.

{\bf Remark:} We observe that, given the choice of $r$ in (\ref{eqn:define_r}), the upper bound on $\mu(\bm{\Psi})$ in (\ref{eqn:coherence_bound}) decays like $1/\sqrt{N}$ for asymptotically large $N$, which is consistent with the Central Limit Theorem.

In order to prove Theorem \ref{the:mu}, we first need to compute the maximum of the Legendre PC basis functions $\psi_{\bm{\alpha}}(\bm{y})$. This is given in the following lemma.

\begin{lemma}[Maximum of $\psi_{\bm{\alpha}}(\bm{y})$]
\label{lem:diameter}
Let $\{\psi_{\bm{\alpha}}(\bm{y})\}$ be the Legendre polynomial chaos basis of total order $p$ in $d$ i.i.d uniform random variables $\bm{y}$ (as defined in (\ref{eqn:multi-index})) and with cardinality $P$. Then,
\begin{equation}
\label{eqn:sup_norm_psi}
\Vert\psi_{\bm{\alpha}}\Vert_{L^{\infty}(\Gamma)}\leq P^{c_{p,d}},
\end{equation}
with a constant
\begin{equation}
\label{eqn:c_pd}
c_{p,d}:=\frac{\ln 3}{2}\frac{p}{\ln\left(\frac{(p+d)!}{p!d!}\right)}.
\end{equation}
\end{lemma}
{\it Proof.} Given the equality
\begin{equation}
\label{eqn:inf_norm_legendre}
\Vert \psi_{\alpha_i}\Vert_{L^{\infty}(\Gamma_i)}=\sqrt{2\alpha_i+1},\nonumber
\end{equation}
\noindent we have
\begin{equation}
\label{eqn:proof_inf_norm_1-gh}
\Vert \psi_{\bm{\alpha}}\Vert_{L^{\infty}(\Gamma)}=\prod_{i=1}^{d}\sqrt{2\alpha_i+1}.\nonumber
\end{equation}
\noindent Using the constraint $\alpha_i\in \mathbb{N}_{0}^{1}$ and $\sum_{i=1}^d \alpha_i\leq p$, the right-hand-side is maximized when $p$ of $\alpha_i$s are equal to one and $d-p$ equal to zero (we assume $d\geq p$). We deduce that
\begin{equation}
\label{eqn:proof_inf_norm_1}
\Vert \psi_{\bm{\alpha}}\Vert_{L^{\infty}(\Gamma)}\leq 3^{\frac{p}{2}}=P^{\frac{\ln 3}{2} \frac{p}{\ln P}}=P^{c_{p,d}}.\nonumber
\end{equation}

We are now ready to prove Theorem \ref{the:mu}.

{\it Proof of Theorem \ref{the:mu}.} The mutual coherence $\mu(\bm{\Psi})$ is
\begin{equation}
\label{eqn:mutual_coherence_legendre}
\mu(\bm{\Psi})=\max_{1\leq j,k\leq P,\;j\neq k}\frac{\left\vert \frac{1}{N}\sum_{i=1}^{N}\psi_{\bm{\alpha}_j}(\bm{y}_i)\psi_{\bm{\alpha}_k}(\bm{y}_i)\right\vert}{\left(\frac{1}{N}\sum_{i=1}^{N}\psi_{\bm{\alpha}_j}^{2}(\bm{y}_i)\right)^{1/2}\left(\frac{1}{N}\sum_{i=1}^{N}\psi_{\bm{\alpha}_k}^{2}(\bm{y}_i)\right)^{1/2}}.
\end{equation}
\noindent Given the independence of samples $\{\bm{y}_i\}_{i=1}^{N}$ and using the McDiarmid's inequality, we obtain
\begin{equation}
\label{eqn:coherence_1a-1}
Prob\left[ \left\vert \frac{1}{N}\sum_{i=1}^{N}\psi_{\bm{\alpha}_j}(\bm{y}_i)\psi_{\bm{\alpha}_k}(\bm{y}_i)\right\vert \ge r \right]\leq 2\exp\left({\frac{-2Nr^2}{(2\|\psi_{\bm{\alpha}_j}(\bm{y}_i)\psi_{\bm{\alpha}_k}(\bm{y}_i)\|_{L^\infty(\Gamma_i)})^2}}\right).
\end{equation}
Using Lemma \ref{lem:diameter}, we may probabilistically bound the numerator in (\ref{eqn:mutual_coherence_legendre}) as
\begin{equation}
\label{eqn:coherence_1}
Prob\left[ \left\vert \frac{1}{N}\sum_{i=1}^{N}\psi_{\bm{\alpha}_j}(\bm{y}_i)\psi_{\bm{\alpha}_k}(\bm{y}_i)\right\vert \ge r \right]\leq 2\exp\left({\frac{-Nr^2}{2P^{4c_{p,d}}}}\right),\nonumber
\end{equation}
\noindent in which we exploit the orthonormality of $\psi_{\bm{\alpha}_j}(\bm{y})$ and $\psi_{\bm{\alpha}_k}(\bm{y})$, i.e., $\mathbb{E}[\psi_{\bm{\alpha}_j}\psi_{\bm{\alpha}_k}]=\delta_{jk}$. Similarly, for $j=1,\cdots,P$, we have
\begin{equation}
\label{eqn:coherence_2}
Prob\left[ \frac{1}{N}\sum_{i=1}^{N}\psi_{\bm{\alpha}_j}^{2}(\bm{y}_i) \leq 1-r \right]\leq \exp\left({\frac{-Nr^2}{P^{4c_{p,d}}}}\right)\leq \exp\left({\frac{-Nr^2}{2P^{4c_{p,d}}}}\right).\nonumber
\end{equation}
\noindent Therefore,
\begin{equation}
\label{eqn:coherence_3}
Prob\left[ \frac{\left\vert \frac{1}{N}\sum_{i=1}^{N}\psi_{\bm{\alpha}_j}(\bm{y}_i)\psi_{\bm{\alpha}_k}(\bm{y}_i)\right\vert}{\left(\frac{1}{N}\sum_{i=1}^{N}\psi_{\bm{\alpha}_j}^{2}(\bm{y}_i)\right)^{1/2}\left(\frac{1}{N}\sum_{i=1}^{N}\psi_{\bm{\alpha}_k}^{2}(\bm{y}_i)\right)^{1/2}}\geq \frac{r}{1-r}\right]\leq 4\exp\left({\frac{-Nr^2}{2P^{4c_{p,d}}}}\right)
\end{equation}
\noindent and
\begin{equation}
\label{eqn:coherence_4}
Prob\left[ \mu(\bm{\Psi})\geq \frac{r}{1-r}\right]\leq 4P^2\exp\left({\frac{-Nr^2}{2P^{4c_{p,d}}}}\right).
\end{equation}
Taking
\begin{equation}
\label{eqn:select_r}
r=2\sqrt{\zeta P^{4c_{p,d}}(\ln P)/N},
\end{equation}
\noindent for some $\zeta>1$, we arrive at the statement of the Theorem \ref{the:mu}. $\square$

To summarize, we observe that with large probability, the mutual coherence $\mu(\bm \Psi)$ of the measurement matrix $\bm{\Psi}$ in (\ref{eqn:underdetermind_system}) can be arbitrarily bounded from above by increasing the number $N$ of independent random solution samples. Moreover, given the discussions of Section \ref{sec:sparse_solution}, we know that the solution to problem (\ref{eqn:spde}) is sparse in the Legendre PC basis. These are the two key factors affecting the stability and accuracy of our sparse approximation.

Following \citep{Donoho06a,Bruckstein09}, we next state a {\it sufficient} condition on the sparsity of $u(\bm x,\bm y)$ (or, equivalently, the mutual coherence of $\bm\Psi$) such that the problem $(P_{1,\delta})$ leads to a stable and accurate sparse solution. By stability, we simply mean that the PC coefficients $\bm{c}^{1,\delta}$ recovered from the $\ell_1$-minimization problem $(P_{1,\delta})$ do not blow up in the presence of the truncation error $\delta$. We first assume that $\bm x$ is a fixed point in space and subsequently extend the analysis over the entire spatial domain $\mathcal{D}$.

\begin{lemma}[A condition on sparsity for stability of $(P_{1,\delta})$]
\label{lem:stability_discrete}
Let $u_p^{0}(\bm y):=\sum_{\bm\alpha\in\Lambda_{p,d}}c_{\bm\alpha}^{0}\psi_{\bm\alpha}(\bm y)$ with $c_{\bm\alpha}^{0}=0$ for $\bm\alpha\notin\Lambda_{p_\epsilon,\nu_\epsilon}$ be the sparse Legendre PC approximation of $u(\bm x,\bm y)$ at a spatial point $\bm x$ where the sparse index set $\Lambda_{p_\epsilon,\nu_\epsilon}$ is defined in \eqref{eqn:sparse-set-epsilon}. Assume that the vector of PC coefficients $\bm{c}^{0}$ satisfies the sparsity condition
\begin{equation}
\label{eqn:sparsity-condition}
\Vert\bm{c}^{0}\Vert_{0}=\vert\Lambda_{p_\epsilon,\nu_\epsilon}\vert<(1+1/\mu(\bm{\Psi}))/4.
\end{equation}
Let $u_p^{1,\delta}(\bm{y}):=\sum_{\bm \alpha\in\Lambda_{p,d}}c_{\bm\alpha}^{1,\delta}\psi_{\bm\alpha}(\bm y)$ be the approximation of $u(\bm{y})$ with coefficients $\bm{c}^{1,\delta}$ computed from the $\ell_1$-minimization problem $(P_{1,\delta})$.Then, with probability at least $1-\exp(-\frac{N}{8 P^{4 c_{p,d}}})$ (on the $N$ samples $\{u(\bm{y}_i)\}_{i=1}^{N}$) and for all $\delta\geq 0$, the solution $u_p^{1,\delta}$ must obey
\begin{equation}
\label{eqn:stability0}
\mathbb{E}\left[\left(u_p^{0}-u_p^{1,\delta}\right)^2\right] \leq \frac{4}{3N} \frac{\left(\delta+\|\bm{\Psi}\bm{c}^0-\bm{u}\|_2 \right)^2}{1-\mu(\bm{\Psi})\left(4\Vert\bm{c}^0\Vert_{0}-1\right)}.
\end{equation}
\end{lemma}
{\it Proof.} Using Theorem 3.1 of \citep{Donoho06a}, we obtain that if $\bm{c}^0$ satisfies the sparsity condition $\Vert\bm{c}^0\Vert_{0}<(1+1/\mu(\bm{\Psi}))/4$, then
\begin{equation}
\sum_{\bm\alpha} \left(c_{\bm\alpha}^{1,\delta}-c_{\bm\alpha}^{0}\right)^2
\left \Vert\bm{\psi}_{\bm\alpha}\right\Vert_{2}^2 \leq \frac{\left(\delta+\|\bm{\Psi}\bm{c}^0-\bm{u}\|_2 \right)^2}{1-\mu(\bm{\Psi})\left(4\Vert\bm{c}^0\Vert_{0}-1\right)}\nonumber
\end{equation}
where $\left \Vert\bm{\psi}_{\bm\alpha}\right\Vert_{2}$ is the $\ell_2$-norm of the column of $\bm{\Psi}$ corresponding to index $\bm{\alpha}$. The presence of $\left \Vert\bm{\psi}_{\bm\alpha}\right\Vert_{2}$ is due to the fact that the columns of $\bm{\Psi}$ are not normalized.
Next, using McDiarmid's inequality and the independence of the entries $\bm{\Psi}[i,j]$ for distinct $i$s, we obtain that
\begin{equation}
\label{eqn:stability1}
Prob\left[\sum_{\bm\alpha} \left(c_{\bm\alpha}^{1,\delta}-c_{\bm\alpha}^{0}\right)^2
\left \Vert\bm{\psi}_{\bm\alpha}\right\Vert_{2}^2\leq  \frac{3N}{4} \sum_{\bm\alpha} \left(c_{\bm\alpha}^{1,\delta}-c_{\bm\alpha}^{0}\right)^2
\left \Vert\bm{\psi}_{\bm\alpha}\right\Vert_{2}^2\right]
\leq \exp\left(-\frac{N}{8 \|\psi_{\bm{\alpha}}\|_{L^\infty(\Gamma)}^4}\right)
\end{equation}
We conclude using Lemma \ref{lem:diameter} and the fact that due to the orthonormality of $\{\psi_{\bm{\alpha}}(\bm{y})\}$ we have $\mathbb{E}\left[\left(u_p^{0}-u_p^{1,\delta}\right)^2\right] =\left\Vert \bm{c}^0-\bm{c}^{1,\delta} \right\Vert_{2}^{2}$. $\square$

{\bf Remark:} The error bound in \eqref{eqn:stability0} is not tight; in fact, the actual error is significantly smaller than the upper bound given in \eqref{eqn:stability}. More importantly, according to \citep{Donoho06a}, the sparsity condition \eqref{eqn:sparsity-condition} is unnecessarily too restrictive. In practice, both far milder sparsity conditions are needed and much better actual errors are achieved.

{\bf Remark:} We will later use the sparsity condition \eqref{eqn:sparsity-condition} to derive the sufficient condition \eqref{eqn:maximum_sparsity} (together with \eqref{eqn:sparsity_condition}) on the number $N$ of random samples needed for a successful recovery. As the condition \eqref{eqn:sparsity-condition} is too restrictive, the theoretical lower bound on $N$ given in \eqref{eqn:maximum_sparsity} and \eqref{eqn:sparsity_condition} is too pessimistic.

{\bf Remark:} According to Lemma \ref{lem:stability_discrete}, we do not need to know {\it a priori} the sparse index set $\Lambda_{p_\epsilon,\nu_\epsilon}$; only the sparsity condition \eqref{eqn:sparsity-condition} is required.

{\bf Remark:} We stated Lemma \ref{lem:stability_discrete} for the case where the sparsity of the PC expansion is due to the fact that the effective dimensionality $\nu_\epsilon$ is potentially smaller than $d$. However, as far as the stability condition \eqref{eqn:sparsity-condition} is satisfied, similar stability results are valid for situations where dominant basis are defined over all the dimensions.

{\bf Remark:} With slight modifications, a similar argument as in Lemma \ref{lem:stability_discrete} may be asserted for the solution of $\ell_0$-minimization problem $(P_{0,\delta})$. Specifically, in that case, we only require a sparsity limit $\Vert\bm{c}^{0}\Vert_{0}=\vert\Lambda_{p_\epsilon,\nu_\epsilon}\vert<(1+1/\mu(\bm{\Psi}))/2$ to achieve the error estimate $\mathbb{E}\left[\left(u_p^{0}-u_p^{0,\delta}\right)^2\right] \leq \frac{4}{3N} \frac{\left(\delta+\|\bm{\Psi}\bm{c}^0-\bm{u}\|_2 \right)^2}{1-\mu(\bm{\Psi})\left(2\Vert\bm{c}^0\Vert_{0}-1\right)}$.

Notice that the normalized truncation error
\begin{equation}
\label{eqn:truncation_error}
\epsilon_{N}^{2}:=\frac{\|\bm{\Psi}\bm{c}^0-\bm{u}\|_2^2}{N}
\end{equation}
is the sample average estimate of the mean-squares sparse approximation error $\mathbb{E}\left[(u-u_{p}^{0})^{2}\right]$ at the point $\bm x$ and is a function of samples $\{\bm y\}_{i=1}^{N}$ in addition to the order $p_\epsilon$ and the dimensionality $\nu_\epsilon$ of the sparse PC expansion. For a given set of $N$ random independent samples $\{\bm y\}_{i=1}^{N}$ and $\{u(\bm y)\}_{i=1}^{N}$, we may  bound $\epsilon_N^{2}$ in probability using McDiarmid's inequality, i.e.,
\begin{equation}
\label{eqn:bound_truncation_error}
Prob\left[\epsilon_N^2 \geq \mathbb{E}\left[(u-u_{p}^{0})^{2}\right] + r\right]\leq \exp\left(-2N \frac{r^2}{\|u - u_{p}^0\|_{L^\infty(\Gamma)}^4} \right).
\end{equation}
Although our sparse approximations are point-wise in space, we are ultimately interested in deriving suitable global stability and error estimates for our sparse reconstructions. Such extensions are readily available from Lemma \ref{lem:stability_discrete} and are stated in the following corollary.

\begin{corollary}
\label{cor:stability_continuous}
Let $u_p^{0}(\bm x, \bm y):=\sum_{\bm\alpha\in\Lambda_{p,d}}c_{\bm\alpha}^{0}(\bm x)\psi_{\bm\alpha}(\bm y)$ with $c_{\bm\alpha}^{0}(\bm x)=0$ for $\bm\alpha\notin\Lambda_{p_\epsilon,\nu_\epsilon}$ be the sparse Legendre PC approximation of $u(\bm x,\bm y)$ where the sparse index set $\Lambda_{p_\epsilon,\nu_\epsilon}$ is defined in \eqref{eqn:sparse-set-epsilon}. Assume that the vector of PC coefficients $\bm{c}^{0}(\bm x)$ satisfies the sparsity condition $\Vert\bm{c}^{0}\Vert_{0}=\vert\Lambda_{p_\epsilon,\nu_\epsilon}\vert<(1+1/\mu(\bm{\Psi}))/4$. Let $u_p^{1,\delta}(\bm x,\bm{y}):=\sum_{\bm \alpha\in\Lambda_{p,d}}c_{\bm\alpha}^{1,\delta}(\bm x)\psi_{\bm\alpha}(\bm y)$ be the approximation of $u(\bm x,\bm{y})$ with coefficients $\bm{c}^{1,\delta}(\bm x)$ computed from the $\ell_1$-minimization problem $(P_{1,\delta})$ at any point $\bm x\in \mathcal{D}$.Then, with probability at least $1-\exp(-\frac{N}{8 P^{4 c_{p,d}}})$ (on the $N$ samples $\{u(\bm{y}_i)\}_{i=1}^{N}$) and for all $\delta\geq 0$, the solution $u_p^{1,\delta}$ must obey
\begin{equation}
\label{eqn:stability}
\left\Vert u_p^{0}-u_p^{1,\delta}\right\Vert_{L^2\left(\mathcal{D},L^2(\Gamma)\right)}^2 \leq \frac{4}{3N} \frac{\left(\delta+\left\Vert\bm{\Psi}\bm{c}^0-\bm{u}\right\Vert_{L^{2}(\mathcal{D},\ell_2(\mathbb{R}^{N}))} \right)^2}{1-\mu(\bm{\Psi})\left(4\Vert\bm{c}^0\Vert_{0}-1\right)}.
\end{equation}
Furthermore, the normalized error $\epsilon_{N}^2=:\frac{\left\Vert\bm{\Psi}\bm{c}^0-\bm{u}\right\Vert_{L^{2}\left(\mathcal{D},\ell_2(\mathbb{R}^{N})\right)}^2}{N}$ is  bounded from above in probability  through
\begin{equation}
\label{eqn:bound_truncation_error_global}
Prob\left[\epsilon_N^2 \geq \left\Vert u- u_p^{0}\right\Vert_{L^2\left(\mathcal{D},L^2(\Gamma)\right)}^2 + r\right]\leq\exp\left(-2N \frac{r^2}{\|u - u_{p}^0\|_{L^2\left(\mathcal{D},L^\infty(\Gamma)\right)}^4} \right).
\end{equation}
\end{corollary}

We have now all the necessary tools to proceed with the proof of our main result stated in Theorem \ref{the:main} which is primarily a direct consequence of Lemma \ref{lem:sparsity}, Theorem \ref{the:mu}, and Corollary \ref{cor:stability_continuous}.

\emph{Proof of Theorem \ref{the:main}.} We first note that, given the conditions of Lemma \ref{lem:sparsity}, the solution to problem \eqref{eqn:spde} admits a sparse Legendre PC expansion $u_{p}$ with sparsity $S=\vert \Lambda_{p_\epsilon,\nu_\epsilon}\vert\lesssim\epsilon^{-1/\rho}$ when $\epsilon\geq \exp\left(-\left(\frac{d}{c_1}\right)^\kappa\right)$ for some constants $c_1$ and (arbitrary) $\rho>0$. Notice that the sparse approximation $u_{p}$ has an accuracy better than $\epsilon$ in the $H_{0}^{1}(\mathcal{D},L^{\infty}(\Gamma))$ sense. Based on Theorem \ref{the:mu}, it is sufficient to have random solution samples of size $N\geq64 P^{4c_{p,d}}(\ln P)S$, to meet the sparsity requirement $S=\vert\Lambda_{p_\epsilon,\nu_\epsilon}\vert<(1+1/\mu(\bm{\Psi}))/4$, in Corollary \ref{cor:stability_continuous}, with probability at least $1-4P^{2-2S_{\max}}$ where $S_{\max}:=\frac{N}{64 P^{4c_{p,d}}(\ln P)}$. On the other hand, given $N$ random samples of solution, we require $\epsilon\geq \frac{1}{S_{\max}^\rho}$ to satisfy the sparsity condition. Given Corollary \ref{cor:stability_continuous} and using the triangular and Poincar\'e inequalities, with probability at least $1-P^{-8S_{\max}}$, we have
\vspace{-.5cm}
\begin{eqnarray}
\label{eqn:tri_poin}
\left\Vert u-u_{p}^{1,\delta}\right\Vert_{L^2(\mathcal{D}, L^2(\Gamma))}&&\leq c_{\mathcal{D}}\left\Vert u-u_{p}^{0}\right\Vert_{H_0^{1}(\mathcal{D}, L^2(\Gamma))}+\left\Vert u_{p}^{0}-u_{p}^{1,\delta}\right\Vert_{L^2(\mathcal{D}, L^2(\Gamma))}\nonumber  \\
&&\leq c_{\mathcal{D}}\epsilon+\frac{2}{\sqrt{3N}} \frac{\delta+\left\Vert\bm{\Psi}\bm{c}^0-\bm{u}\right\Vert_{L^{2}(\mathcal{D},\ell_2(\mathbb{R}^{N}))}}{\sqrt{1-\mu(\bm{\Psi})\left(4S-1\right)}},\nonumber
\end{eqnarray}
for all $\delta>0$. Moreover, by choosing $r= \frac{1}{4}\left\Vert u- u_p^{0}\right\Vert_{L^2\left(\mathcal{D},L^\infty(\Gamma)\right)}^2$ in \eqref{eqn:bound_truncation_error_global},
\begin{equation}
\label{eqn:epsilon_n_bound}
\epsilon_{N}^2=:\frac{\left\Vert\bm{\Psi}\bm{c}^0-\bm{u}\right\Vert_{L^{2}\left(\mathcal{D},\ell_2(\mathbb{R}^{N})\right)}^2}{N}\leq \left\Vert u- u_p^{0}\right\Vert_{L^2\left(\mathcal{D},L^2(\Gamma)\right)}^2+ \frac{1}{4}\left\Vert u- u_p^{0}\right\Vert_{L^2\left(\mathcal{D},L^\infty(\Gamma)\right)}^2\leq \frac{5}{4}c_{\mathcal{D}}^2\epsilon^2\nonumber
\end{equation}
with probability at least $1-P^{-8S_{\max}P^{4c_{p,d}}}$. Finally, by taking
\begin{equation}
\label{eqn:c1_c2}
c_4:=c_{\mathcal{D}}\left(1+\frac{\sqrt{5}}{\sqrt{3}\sqrt{1-\mu(\bm{\Psi})\left(4S-1\right)}}\right)\quad \mathrm{and}\quad c_5:=\frac{2}{\sqrt{3}\sqrt{1-\mu(\bm{\Psi})\left(4S-1\right)}},\nonumber
\end{equation}
we arrive at the statement of the Theorem \ref{the:main}. $\square$

\subsection{Choosing the truncation error tolerance $\delta$}
\label{sec:delta}

An important component of the sparse approximation using $(P_{1,\delta})$ and $(P_{0,\delta})$ is the selection of the truncation error tolerance $\delta$. Although the stability bounds given in Lemma \ref{lem:stability_discrete} and Corollary \ref{cor:stability_continuous} are valid for any $\delta\geq0$, the actual error and the sparsity level of the solution to $(P_{1,\delta})$ and $(P_{0,\delta})$ depend on the choice of $\delta$. Ideally, we desire to choose $\delta\approx\Vert\bm\Psi\bm{c}^{0}-\bm u\Vert_2$; while larger values of $\delta$ deteriorate the accuracy of the approximation, as in Lemma \ref{lem:stability_discrete}, smaller choices of $\delta$ may result in over-fitting the solution samples and, thus, less sparse solutions. In practice, as the exact values of the PC coefficients $\bm{c}^{0}$ are not known, the exact values of the truncation error $\Vert\bm\Psi\bm{c}^{0}-\bm u\Vert$ and, consequently, $\delta$ are not known {\it a priori}. Therefore, $\delta$ has to be estimated, for instance, using statistical techniques such as the cross-validation \citep{Boufounos07,Ward09}.

In this work, we propose a heuristic cross-validation algorithm to estimate $\delta$. We first divide the $N$ available solution samples to $N_r$ reconstruction and $N_v$ validation samples such that $N=N_r+N_v$. The idea is to repeat the solution of $(P_{1,\delta})$ (or $(P_{0,\delta})$) on the reconstruction samples and with multiple values of truncation error tolerance $\delta_r$. We then set $\delta=\sqrt{\frac{N}{N_r}}\hat{\delta}_r$ in which $\hat{\delta}_r$ is such that the corresponding truncation error on the $N_v$ validation samples is minimum. This is simply motivated by the fact that the truncation error on the validation samples is large for values of $\delta_r$ considerably larger and smaller than $\Vert\bm\Psi\bm{c}^{0}-\bm u\Vert_2$ evaluated using the reconstruction samples. While the former is expected from the upper bound on the approximation error in Lemma \ref{lem:stability_discrete}, the latter is due to the over-fitting the reconstruction samples. The following exhibit outlines the estimation of $\delta$ using the above cross-validation approach:\\

\fbox{
\parbox{13cm}{
{\it Algorithm for cross-validation estimation of $\delta$:}\\

\begin{itemize}
\item Divide the $N$ solution samples to $N_r$ {\it reconstruction} and $N_v$ {\it validation} samples. \\[-.3cm]
\item Choose multiple values for $\delta_r$ such that the exact truncation error $\Vert\bm\Psi\bm{c}^{0}-\bm u\Vert_2$ of the reconstruction samples is within the range of $\delta_r$ values. \\[-.3cm]
\item For each value of $\delta_r$, solve $(P_{1,\delta})$ (or $(P_{0,\delta})$) on the $N_r$ reconstruction samples. \\[-.3cm]
\item For each value of $\delta_r$, compute the truncation error $\delta_v:=\Vert\bm\Psi\bm{c}^{1,\delta_r}- \bm u\Vert_2$ (or $\delta_v:=\Vert\bm\Psi\bm{c}^{0,\delta_r}-\bm u\Vert_2$) of the $N_v$ validation samples. \\[-.3cm]
\item Find the minimum value of $\delta_v$ and its corresponding $\hat{\delta}_r:=\delta_r$.\\[-.3cm]
\item Set $\delta=\sqrt{\frac{N}{N_r}}\hat{\delta}_r$.
\end{itemize}
}
}

In the numerical experiments of Section \ref{sec:examples}, we repeat the above cross-validation algorithm for multiple replications of the reconstruction and validation samples. The estimate of $\delta=\sqrt{\frac{N}{N_r}}\hat{\delta}_r$ is then based on the value of $\hat{\delta}_r$ for which the average of the corresponding truncation errors $\delta_v$, over all replications of the validation samples, is minimum. This resulted in more accurate solutions in our numerical experiments.

\subsection{Algorithms}
\label{sec:solvers}
There are several numerical algorithms for solving problems $(P_{0,\delta})$ and $(P_{1,\delta})$ each with different optimization kernel, computational complexity, and degree of accuracy. An in-depth discussion on the performance of these algorithms is outside the scope of the present work; however, below we name some of the available options for each problem and briefly describe the algorithms that have been utilized in our numerical experiments. For comprehensive discussions on this subject, the interested reader is referred to \citep{Bruckstein09,Berg08,Figueiredo07,Becker09,Yang09,Tropp10a}.

{\it Problem $(P_{0,\delta})$:} A brute force search through all possible support sets in order to identify the correct sparsity for the solution $\bm{c}^{0,\delta}$ of $(P_{0,\delta})$ is NP-hard and not practical. Greedy pursuit algorithms form a major class of schemes to tackle the solution of $(P_{0,\delta})$ with a tractable computational cost. Instead of performing an exhaustive search for the support of the sparse solution, these solvers successively find one or more components of the solution that result in the largest improvement in the approximation. Some of the standard greedy pursuit algorithms are Orthogonal Marching Pursuit (OMP) \citep{Pati93,Davis97}, Regularized OMP (ROMP)\citep{Needell07}, Stagewise OMP (StOMP) \citep{Donoho06c}, Compressive Sampling MP (CoSaMP) \citep{Needell08a}, Subspace Pursuit \citep{Dai09}, and Iterative Hard Thresholding (IHT) \citep{Blumensath09}. Under well-defined conditions, all of the above schemes provide stable and accurate solutions to $(P_{0,\delta})$ in a reasonable time.

In the present study, we employ the OMP algorithm to approximate the solution of $(P_{0,\delta})$. Starting from $\bm{c}^{0,\delta,(0)}=\bm{0}$ and an empty active column set of $\bm\Psi$, at any iteration $k$, OMP identifies only one column to be added to the active column set. The column is chosen such that the $\ell_2$-norm of the residual, $\Vert\bm\Psi\bm{c}^{0,\delta,(k)}-\bm u\Vert_2$, is maximally reduced. Having specified the active column set, a least-squares problem is solved to compute the solution $\bm{c}^{0,\delta,(k)}$. The iterations are continued until the error truncation tolerance $\delta$ is achieved. In general, the complexity of the OMP algorithm is $\mathcal{O}(S\cdot N\cdot P)$ where $S:=\Vert\bm{c}^{0,\delta}\Vert_{0}$ is number of non-zero (dominant) entries of $\bm{c}^{0,\delta}$. The following exhibit depicts an step-by-step implementation of the OMP algorithm. 

\fbox{
\parbox{13cm}{
{\it Orthogonal Matching Pursuit (OMP) Algorithm:}\\

\begin{itemize}
\item Set $k=0$.
	\begin{itemize}	
		\item {\footnotesize Set the initial solution $\bm{c}^{0,\delta,(0)}=\bm{0}$ and residual $\bm{r}^{(0)}=\bm{u}-\bm{\Psi}\bm{c}^{0,\delta,(0)}=\bm{u}$}.
		{\footnotesize\item Set the solution support index set $\mathcal{I}^{(0)}=\emptyset$}. \\[-.3cm]
          \end{itemize}
\item  While $\Vert \bm{u}-\bm{\Psi c}^{0,\delta,(k)}\Vert_{2}>\delta$ perform:
	 \begin{itemize}
	 	\item For all $j\notin\mathcal{I}^{(k)}$ evaluate $\epsilon(j)=\Vert \bm{\psi}_j\alpha_j-\bm{r}^{(k)}\Vert_{2}$ with $\alpha_j=\bm{\psi}_j^{T}\bm{r}^{(k)}/\Vert \bm{\psi}_j\Vert_2^{2}$.
		\item Set $k=k+1.$ 
		\item Update the support index set $\mathcal{I}^{(k)}=\mathcal{I}^{(k-1)}\bigcup\left\{\arg\min_{j}\epsilon(j)\right\}$.
		\item Solve for $\bm{c}^{0,\delta,(k)}=\arg\min_{{\bm c}^{0,\delta}}\Vert \bm{u}-\bm{\Psi c}^{0,\delta}\Vert_{2}$ subject to $Support\{\bm{c}^{0,\delta}\}=\mathcal{I}^{(k)}$.
		\item Update the residual $\bm{r}^{(k)}=\bm{u}-\bm{\Psi}\bm{c}^{0,\delta,(k)}$ \\[-.3cm]
	 \end{itemize}
\item Output the solution $\bm{c}^{0,\delta}=\bm{c}^{0,\delta,(k)}$.
\end{itemize}
}
}

Although we chose OMP in our analysis, we note that further studies are needed to identify the most appropriate greedy algorithm for the purpose of this study.

{\it Problem $(P_{1,\delta})$:} The majority of available solvers for $\ell_1$-minimization are based on alternative formulations of $(P_{1,\delta})$, such as the $\ell_1$-norm regularized least-squares problem
\begin{equation}
\label{eqn:QP}
(QP_{\lambda}):\qquad\min_{\bm{c}} \frac{1}{2}\left\Vert\bm{\Psi}\bm{c}-\bm{u}\right\Vert_2^2+\lambda\Vert\bm{Wc}\Vert_{1},
\end{equation}
or the LASSO problem, \citep{Tibshirani94},
\begin{equation}
\label{eqn:LASSO}
(LS_{\tau}):\qquad\min_{\bm{c}} \frac{1}{2}\left\Vert\bm{\Psi}\bm{c}-\bm{u}\right\Vert_2^2\quad\mathrm{subject \;to}\quad\Vert\bm{Wc}\Vert_{1}\leq\tau.
\end{equation}
It can be shown that for an appropriate choice of scalars $\delta$, $\lambda$, and $\tau$, the problems $(P_{1,\delta})$, $(QP_{\lambda})$, and $(LS_{\tau})$ share the same solution \citep{Berg08,Bruckstein09,Tropp10a}.  Among others, the problem $(QP_{\lambda})$ is of particular interest as it is an unconstraint optimization problem. Numerous solvers based on the {\it active set} \citep{Osborne00,Efron04}, {\it interior-point continuation} \citep{Chen01a,Kim07a} and {\it projected gradient} \citep{Daubechies04,Combettes05,Hale08,Bredies08,Bioucas-Dias07,Berg08,Beck09,Becker09} methods have been developed for solving the above formulations of the $\ell_1$-minimization problem.

In our numerical experiments, we adopt the Spectral Projected Gradient algorithm (SPGL1) proposed in \citep{Berg08} and implemented in the MATLAB package {\tt SPGL1} \citep{spgl1:2007} to solve the $\ell_1$-minimization problem $(P_{1,\delta})$ in \eqref{eqn:P1_delta}. SPGL1 is based on exploring the so-called Pareto curve, describing the tradeoff between the $\ell_2$-norm of the truncation error $\Vert\bm{\Psi}\bm{c}-\bm u\Vert_{2}$ and the (weighted) $\ell_1$-norm of the solution $\Vert\bm{Wc}\Vert_{1}$, for successive solution iterations. At each iteration, the LASSO problem \eqref{eqn:LASSO} is solved using the spectral projected gradient technique with a worst-case complexity of $\mathcal{O}(P\ln P)$ where $P$ is the number of columns in $\bm{\Psi}$. Given the error tolerance $\delta$, a scalar equation is solved to identify a value for $\tau$ such that the $(LS_{\tau})$ solution of \eqref{eqn:LASSO} is identical to that of $(P_{1,\delta})$ in \eqref{eqn:P1_delta}. Besides being efficient for large-scale systems where $\bm{\Psi}$ may not be available explicitly, the SPGL1 algorithm is specifically effective for our application of interest as the truncation error $\Vert\bm{\Psi}\bm{c}-\bm u\Vert_{2}$ is known only approximately.

In the next section, we explore some aspects of the proposed scheme through its application to a $1$-$D$ (in space) elliptic stochastic PDE with high-dimensional random diffusion coefficients.

\section{Numerical examples}
\label{sec:examples}

We consider the solution of a one-dimensional, i.e., $D=1$, version of problem (\ref{eqn:spde}),
\begin{eqnarray}
\label{eqn:1D_spde}
&&-\frac{d}{dx}\left(a(x,\omega) \frac{du(x,\omega)}{dx}\right)=1,\quad x\in \mathcal{D}=(0,1),\\
&&\quad u(0,\omega)=u(1,\omega)=0,\nonumber
\end{eqnarray}
where the stochastic diffusion coefficient $a(x,\omega)$ is given by the expansion
\begin{equation}
\label{eqn:example_kle_a}
a(x,\omega)=\bar{a}+\sigma_{a}\sum_{i=1}^{d}\sqrt{\lambda_i}\phi_{i}(x)y_{i}(\omega).
\end{equation}
Here, $\{\lambda_i\}_{i=1}^{d}$ and $\{\phi_i(x)\}_{i=1}^{d}$ are, respectively, $d$ largest eigenvalues and the corresponding eigenfunctions of the Gaussian covariance kernel
\begin{equation}
\label{eqn:gaussian_kernel}
C_{aa}(x_1,x_2)=\exp\left[-\frac{(x_1-x_2)^2}{l_c^2}\right],
\end{equation}
in which $l_c$ is the correlation length of $a(x,\omega)$ that prescribes the decay of the spectrum of $C_{aa}$ in (\ref{eqn:gaussian_kernel}). Random variables $\{y_i(\omega)\}_{i=1}^{d}$ are assumed to be independent and uniformly distributed on $[-1,1]$. The coefficient $\sigma_{a}$ controls the variability of $a(x,\omega)$.

We verify the accuracy and efficiency of the present sparse approximation schemes for both moderate and high-dimensional diffusion coefficient $a(x,\omega)$. These two cases are obtained, respectively, by assuming $(l_c,d) = (1/5,14)$ and $(l_c,d) = (1/14,40)$ in (\ref{eqn:gaussian_kernel}) and (\ref{eqn:example_kle_a}). We further assume that $\bar{a}=0.1$,  $\sigma_a=0.03$ when $d=14$, and $\sigma_a=0.021$ when $d=40$. These choices ensure that all realizations of $a(x,\omega)$ are strictly positive on $\mathcal{D}=(0,1)$. Table \ref{tab:parameters} summarizes the assumed parameters for the two test cases.
\begin{table}[htb]
\begin{center}
\caption{Choices of parameters defining the stochastic description of diffusion coefficient $a(x,\omega)$ in Eq. (\ref{eqn:example_kle_a}).}
 \begin{tabular}{p{1in}p{1in}p{1in}p{1in}p{.5in}}
 \hline
 Case   &    $\bar{a}$    &       $\sigma_{a}$       &   $l_c$  &      $d$       \\
 \hline
     I       &    $0.1$          &                $0.030$       &   $1/5$    &      $14$    \\
     lI      &    $0.1$          &                $0.021$       &   $1/14$  &      $40$     \\
 \hline
 \label{tab:parameters}
 \end{tabular}
 \end{center}
 \end{table}

For both cases, the spatial discretization is done by the Finite Element Method using quadratic elements. A mesh convergence analysis is performed to ensure that spatial discretization errors are inconsequential.

The solution statistics are computed using the conventional Monte Carlo simulation, the isotropic sparse grid stochastic collocation with the Clenshaw-Curtis abscissas \citep{Xiu05a,Babuska07a}, and the proposed sparse approximation techniques. We use the Basis Pursuit Denoising (BPDN) solver implemented in {\tt SPGL1} \citep{spgl1:2007,Berg08} to solve the $\ell_1$-minimization problem $(P_{1,\delta})$ and the Orthogonal Matching Pursuit (OMP) solver in {\tt SparseLab}  \citep{Sparselab} to approximate the $\ell_{0}$-minimization problem $(P_{0,\delta})$. We compare the errors in the mean, standard deviation, and root mean-squares of the solution error at $x=0.5$ using the above methods. The details of the analysis are reported below.

\subsection{Case I: d=14}
\label{sec:case_d_14}

We consider an increasing number $N=\{29,120,200,280,360,421,600\}$ of random solution samples to evaluate the solution $u$ at $x=0.5$ and, consequently, to compute the PC coefficients of the solution using $\ell_1$- and $\ell_0$-minimization. These samples are nested in the sense that we recycle the previous samples when we perform calculations with larger sample sizes. The nested sampling property of our scheme is of paramount importance in large scale calculations where the computational cost of each solution evaluation is enormous. We note that sample sizes $N=29$ and $N=421$, respectively, correspond to the number of nested abscissas in the level $l=1$ and level $l=2$ of the isotropic stochastic collocation with the Clenshaw-Curtis rule.

As elucidated in Section \ref{sec:sucess_rec}, the accuracy of our sparse reconstruction depends on the mutual coherence $\mu(\bm{\Psi})$, the sample size $N$, and the truncation error $\Vert\bm{\Psi c}-\bm{u}\Vert_2$ (hence $\delta$). In order to reduce the approximation error, we need to reduce $\Vert\bm{\Psi c}-\bm{u}\Vert_2$, which may be done by increasing $p$ and, therefore, $P$. However, with a fixed number $N$ of samples, an increase in $P$ may result in a larger mutual coherence and, thus, the degradation of the reconstruction accuracy. Therefore, in practice, we start by approximating the lower order PC expansions when $N$ is small and increase $p$ when larger number of samples become available. Notice that such an adaptivity with respect to the order $p$ is a natural way of refining the accuracy of PC expansions, for instance, when the intrusive stochastic Galerkin scheme is adopted \citep{Ghanem03}. In particular, in this example, for sample sizes $N=\{29,120\}$, we attempt to estimate the coefficients of the $3$rd-order Legendre PC expansion, i.e. $p=3$ and $P=680$. For larger sample sizes $N$, we also include the first $320$ basis function from the $4$th-order chaos, thus resulting in $P=1000$. {\color{black}Since all of the $4$th-order basis functions are not employed, we need to describe the ordering of our basis construction}. We sort the elements of $\{\psi_{\bm{\alpha}}(\bm{y})\}$ such that, for any given order $p$, the random variables $y_i$ with smaller indices $i$ contribute first in the basis.

For each analysis, we estimate the truncation error tolerance $\delta$ based on the cross-validation algorithm described in Section \ref{sec:delta}. For each $N$, we use $N_r\approx 3N/4$ of the samples (reconstruction set) to compute the PC coefficients $\bm{c}^{1,\delta_r}$ and the rest of the samples (validation set) are used to evaluate the truncation error $\delta_{v}$. The cross-validation is performed for four replications of reconstruction and validation sample sets. We then find the value $\hat{\delta}_{r}$ that minimizes the average of $\delta_{v}$ over the four replications of the cross-validation samples. Given an estimate of the truncation error tolerance $\delta\approx\sqrt{4/3}\hat{\delta}_{r}$, we then use all $N$ samples to compute the coefficients $\bm{c}^{1,\delta}$.

\begin{figure}[h]
    \centering
    \begin{tabular}{c}

      \psfrag{x}{$x$}	
      \includegraphics[width=5.0in]{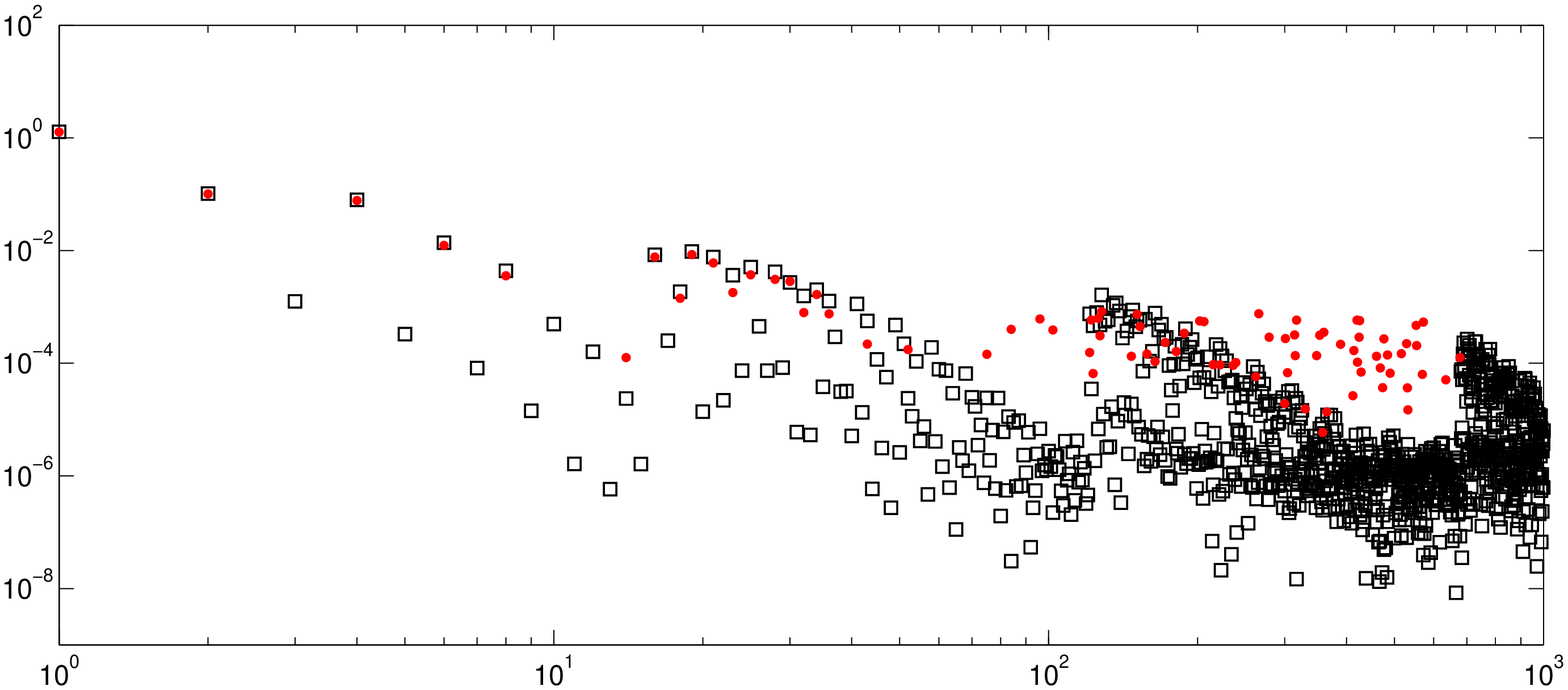}
      \put(-345,60){\begin {sideways} $\vert c_{\bm \alpha_i}\vert$ \end{sideways}}
      \put(-300,25){$\bm{(a)}$}
      \\
      \psfrag{x}{$x$}	
      \includegraphics[width=5.0in]{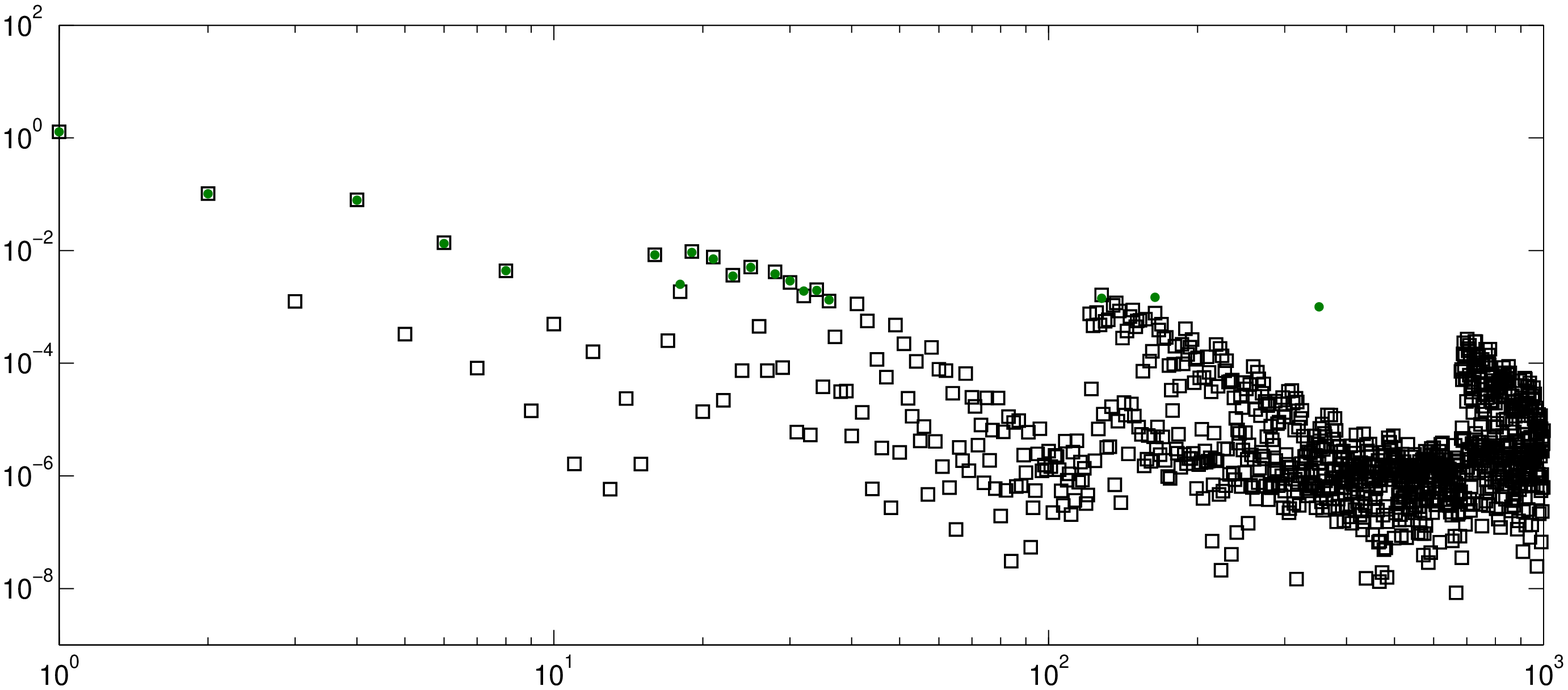}
      \put(-345,60){\begin {sideways} $\vert c_{\bm \alpha_i}\vert$ \end{sideways}}
      \put(-300,25){$\bm{(b)}$}
      \\
      \psfrag{x}{$x$}
      \includegraphics[width=5.0in]{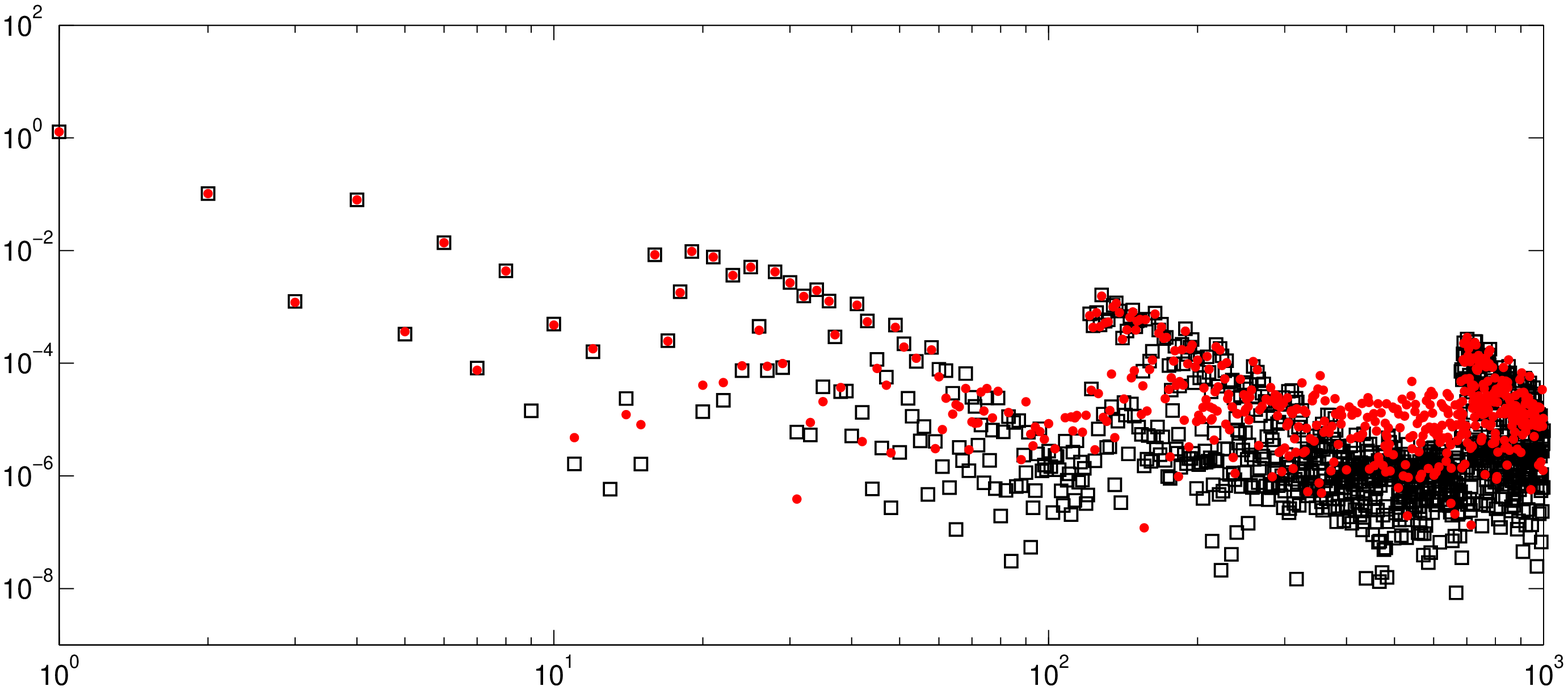}
      \put(-345,60){\begin {sideways} $\vert c_{\bm \alpha_i}\vert$ \end{sideways}}
      \put(-300,25){$\bm{(c)}$}
      \\
      \psfrag{x}{$x$}
      \includegraphics[width=5.0in]{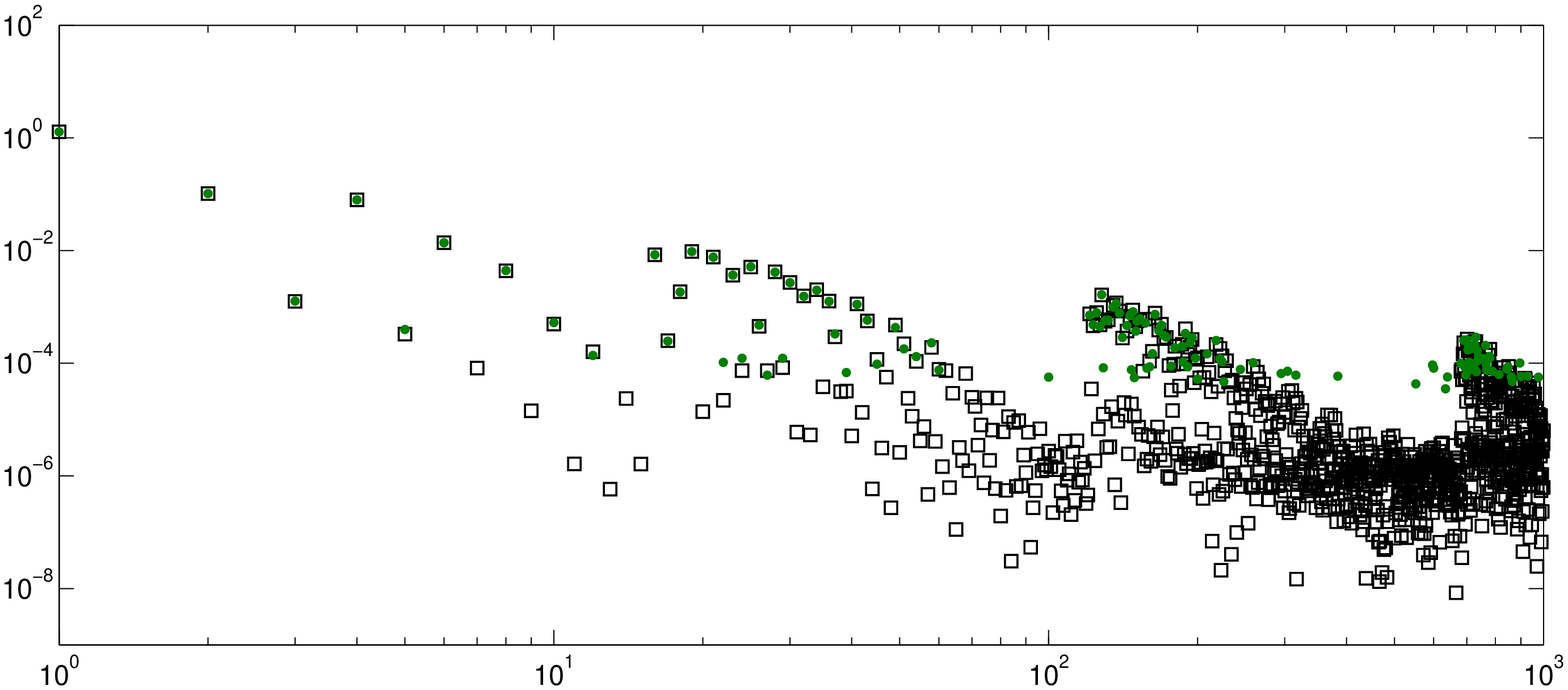}
      \put(-345,60){\begin {sideways} $\vert c_{\bm \alpha_i}\vert$ \end{sideways}}
      \put(-230,-4){\footnotesize Index of PC coefficients ($i$)}
      \put(-300,25){$\bm{(d)}$}
     \end{tabular}
          \caption{Approximation of polynomial chaos (PC) coefficients $\bm c$ of $u(0.5,\bm{y})$ using BPDN and OMP for $d=14$. (a) BPDN with $N=120$ samples, (b) OMP with $N=120$ samples, (c) BPDN with $N=600$ samples, and (d) OMP with $N=600$ samples . `Exact' coefficients computed from level 8 stochastic collocation with the Clenshaw-Curtis abscissas ({\scriptsize $\square$}); BPDN and OMP ($\bullet$).}
     \label{fig:pce_coef_14}
\end{figure}

Figure \ref{fig:pce_coef_14} compares the `exact' PC coefficients with those obtained using BPDN and OMP solvers. We only demonstrate the results corresponding to sample sizes $N=\{120,600\}$. An `exact' solution is computed using the level $l=8$ stochastic collocation for which the approximation errors are negligible in our comparisons. We observe that BPDN tends to give less sparse solutions compared to OMP. This is due to the facts that $i)$ the solution is not exactly sparse, i.e., there are many non-zero (although negligible) coefficients $c_{\bm \alpha}$, $ii)$ the $\ell_1$ cost function does not impose a sufficiently large penalty on the small coefficients as does the $\ell_0$ cost function, and $iii)$ the truncation error tolerance $\delta$ may be under-estimated. To reduce this issue, a number of modifications, including the reweighted $\ell_1$-minimization \citep{Candes08a,Xu09a,Needell09,Khajehnejad10}, have been introduced in the literature that are the subjects of our future work. In contrary, OMP results in more sparse solutions as it adds basis function one-at-a-time until the residual falls below the truncation error. However, as is seen in Figs. \ref{fig:pce_coef_14} (b) and (d), a number of small coefficients are still over-estimated. This is primarily due to under-estimation of the truncation error tolerance $\delta$ in the cross-validation algorithm.

\clearpage

\begin{figure}
    \centering
    \begin{tabular}{cc}

      \psfrag{x}{$x$}	
      \includegraphics[width=2.67in]{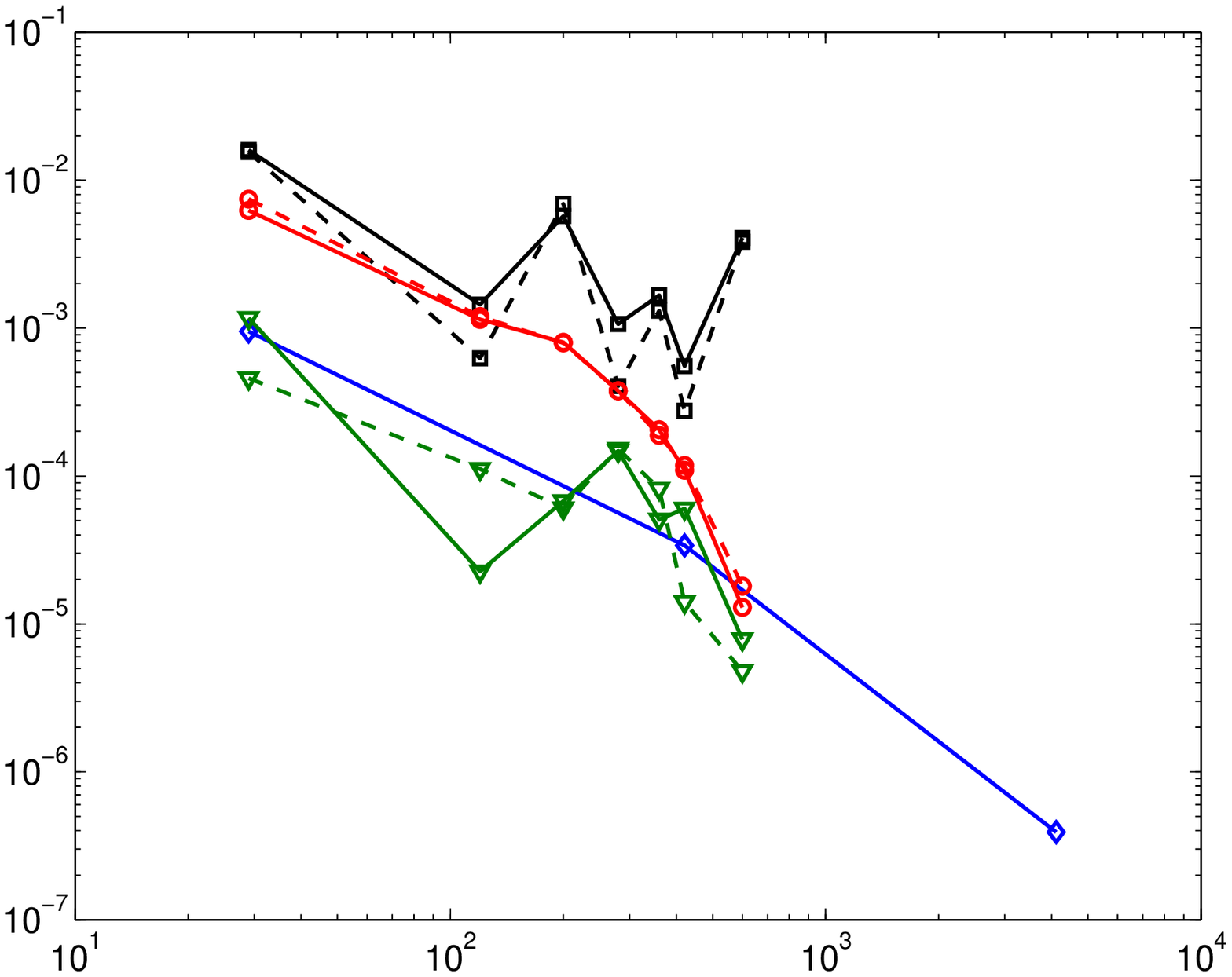}
      \put(-205,35){\begin {sideways} {\footnotesize Rel. error in mean} \end{sideways}}
      \put(-95,-10){$N$}
      \put(-30,130){$\bm{(a)}$}
	&
      \psfrag{x}{$x$}	
      \includegraphics[width=2.67in]{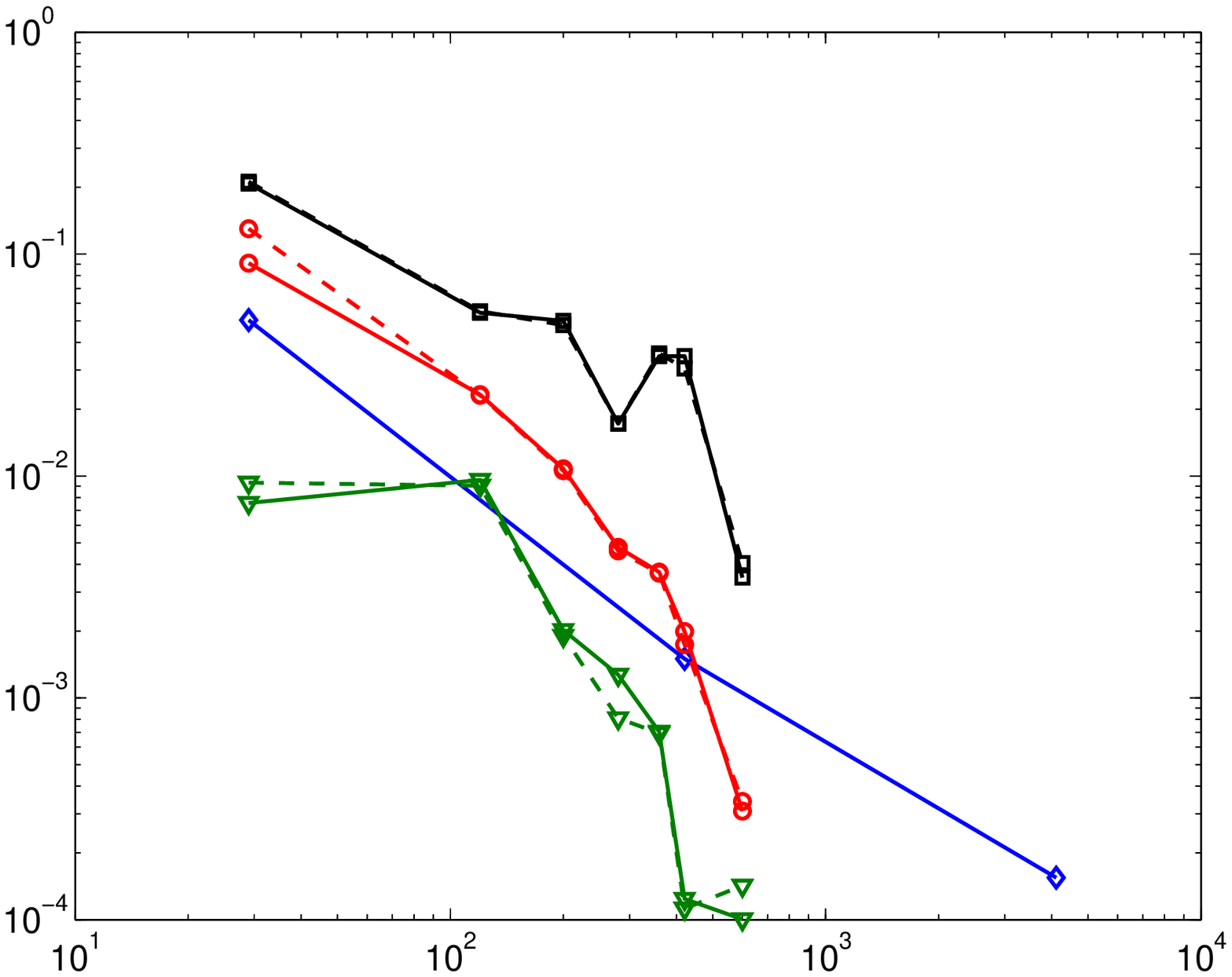}
      \put(-205,35){\begin {sideways} {\footnotesize Rel. error in s.d.} \end{sideways}}
      \put(-95,-10){$N$}
      \put(-30,130){$\bm{(b)}$}
      \\
      \psfrag{x}{$x$}	
      \includegraphics[width=2.67in]{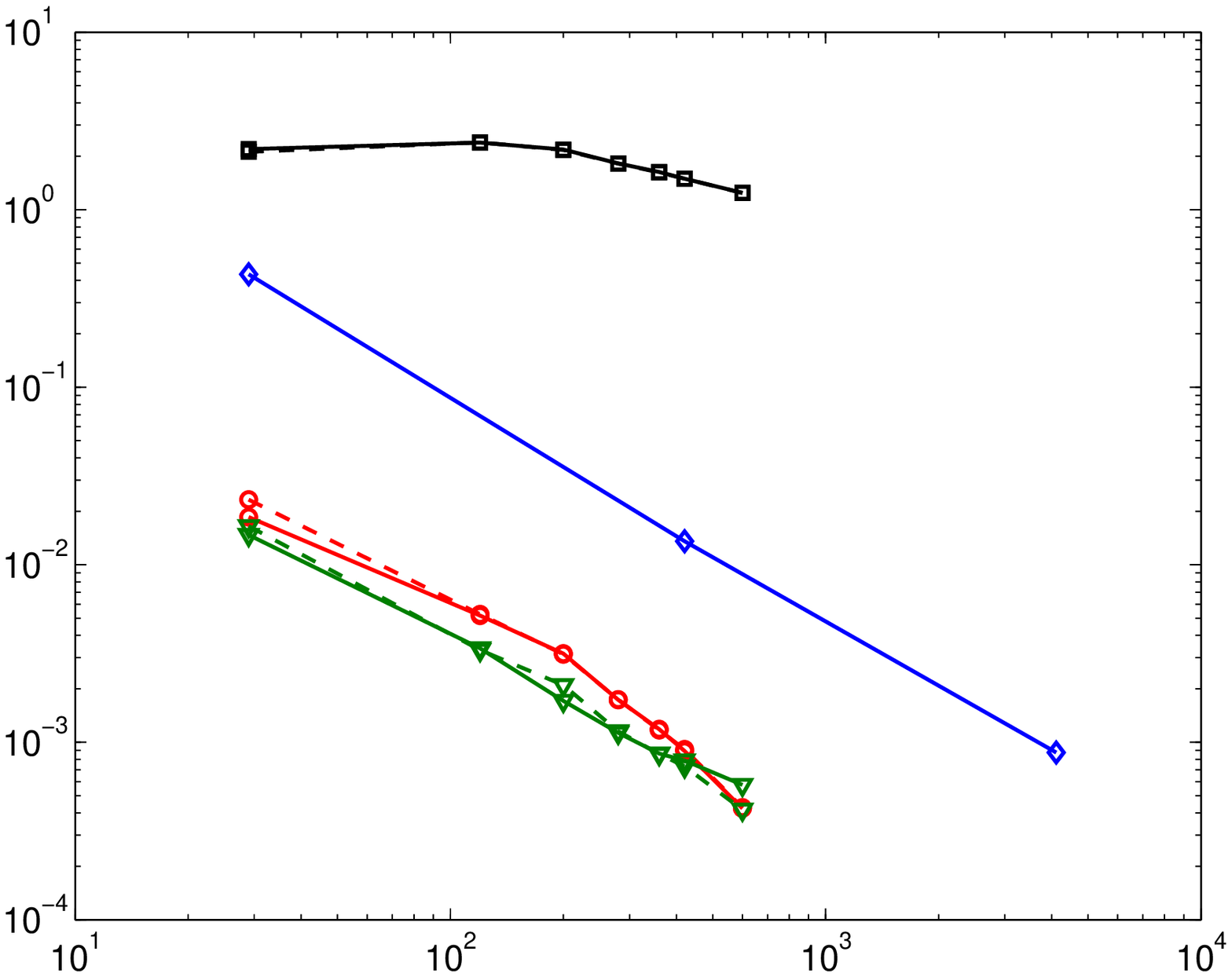}
      \put(-205,40){\begin {sideways} {\footnotesize Rel. rms error} \end{sideways}}
      \put(-95,-10){$N$}
      \put(-30,130){$\bm{(c)}$}
	&
      \psfrag{x}{$x$}	
      \includegraphics[width=2.95in]{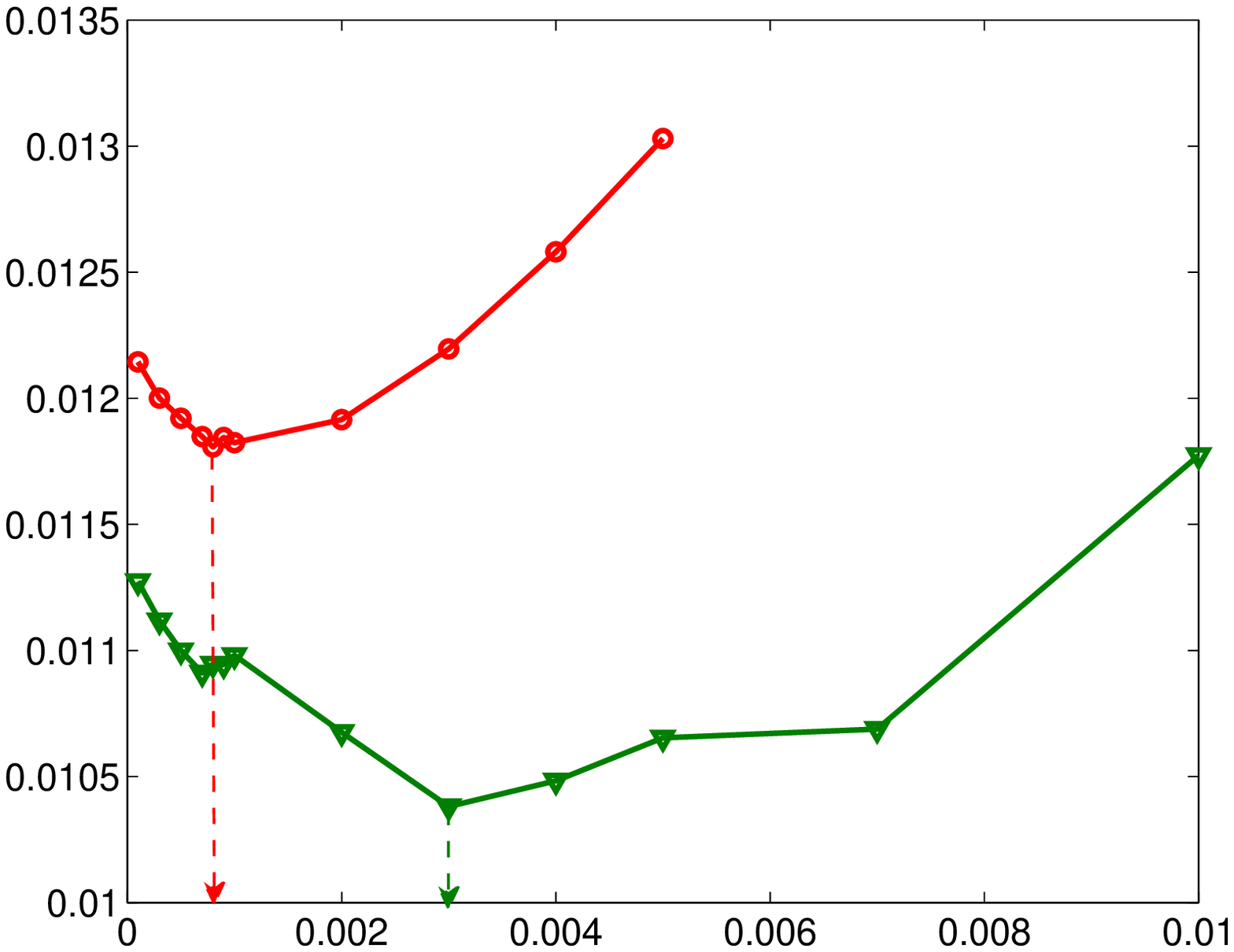}
      \put(-215,75){\begin {sideways} $\delta_{v}$ \end{sideways}}
      \put(-100,0){$\delta_{r}$}
      \put(-170,140){\footnotesize $N=600, p=4$}
      \put(-135,8){\scriptsize \color{OliveGreen} $\hat{\delta}_{r}$}
      \put(-172,8){\scriptsize \color{red} $\hat{\delta}_{r}$}
      \put(-45,130){$\bm{(d)}$}
     \end{tabular}
         \caption{Comparison of relative error in solution statistics at $x=0.5$ for the Monte Carlo simulation, isotropic sparse grid stochastic collocation with the Clenshaw-Curtis abscissas, and the proposed sparse approximations (BPDN and OMP) for $d=14$. Two sets of independent random samples of $u(0.5,\bm{y})$ are generated first and are used for the Monte Carlo simulation, BPDN, and OMP. The solid and dashed lines correspond to the first and second sets of samples, respectively. (a) Relative error in the mean; (b) Relative error in the standard deviation; (c) Relative root-mean-squares (rms) error; and (d) Estimation of $\delta$ using cross-validation: $\delta$ is computed from $\delta_{r}$ for which $\delta_{v}$ is minimum. (Monte Carlo simulation ({\scriptsize $\square$}); stochastic collocation ($\diamond$); BPDN ($\circ$); OMP ($\triangledown$)).}
    \label{fig:mean_sd_mse_cross_14}
\end{figure}

The convergence of the mean, standard deviation, and root mean-squares of the approximation error for  $u(0.5,\bm{y})$ is illustrated in Figs. \ref{fig:mean_sd_mse_cross_14} (a), (b), and (c), respectively. For the case of stochastic collocation, we apply sparse grid quadrature (cubature) integration rule to directly compute the mean and the standard deviation. The root mean-squares error of the Monte Carlo and the stochastic collocation solution are evaluated by estimating the corresponding PC coefficients using sampling and sparse grid quadrature integration, respectively, and then comparing them with the exact coefficients. 

\clearpage

To make a meaningful comparison, for each $N$, the samples used to compute the solution statistics by the conventional Monte Carlo, BPDN,  and OMP are identical. In this sense, the sparse approximation using $\ell_1$- and $\ell_0$-minimizations may be viewed as only post-processing steps in the Monte Carlo simulation. As the sample sizes are finite, the estimates of the PC coefficients $\bm c$ are sample dependent and are in fact random. To demonstrate the convergence of the algorithm with respect to different sets of samples, we repeat the analysis for two independent sets of $N$ samples and report the corresponding statistics errors with solid and dashed lines. Although for different solution samples of size $N$ the estimates of $\bm c$ and the solution statistics are not identical, the approximation converges, with large probability, for any set of samples with sufficiently large size $N$ (see Theorem \ref{the:main}). 

Figure \ref{fig:mean_sd_mse_cross_14} (d) illustrates the statistical estimation of $\delta$ using the cross-validation approach described in Section \ref{sec:delta}. The estimation of  $\delta$ is slightly different in BPDN and OPM, this is a consequence of different reconstruction accuracy of these techniques. Moreover, the solution of the BPDN algorithm is less sensitive to small perturbations in the truncation error $\delta$ compared to that of the OMP algorithm. This is justified by the fact that the $\ell_0$-norm is highly discontinuous. 

{\bf Remark:} Despite the conventional implementation of the stochastic collocation approach where the approximation refinement requires a certain number of extra samples, the $\ell_1$- and $\ell_0$-minimizations may be implemented using arbitrary numbers of additional samples, which is an advantage, particularly, when only a limited number of samples is afforded.

\subsection{Case II: d=40}
\label{sec:case_d_40}

The objective of this example is to highlight that a sparse reconstruction may lead to significant computational savings for problems with high-dimensional random inputs. Similar to the analysis of Case I described in Section \ref{sec:case_d_14}, we compute the solution statistics using multiple numbers of independent samples. More specifically, we evaluate the solution at $x=0.5$ for independent samples of size $N=\{81,200,400,600,800,1000\}$. The number of grid points in the level $l=1$ and $l=2$ of the Clenshaw-Curtis rule in dimension $d=40$ is $N=81$ and $N=3281$, respectively. To obtain a reference solution, the $3$rd order PC coefficients $\bm c$ of the solution at $x=0.5$ are computed using level $l=5$ stochastic collocation with the Clenshaw-Curtis rule.

\begin{figure}
    \centering
    \begin{tabular}{c}

      \psfrag{x}{$x$}
      \includegraphics[width=4.8in]{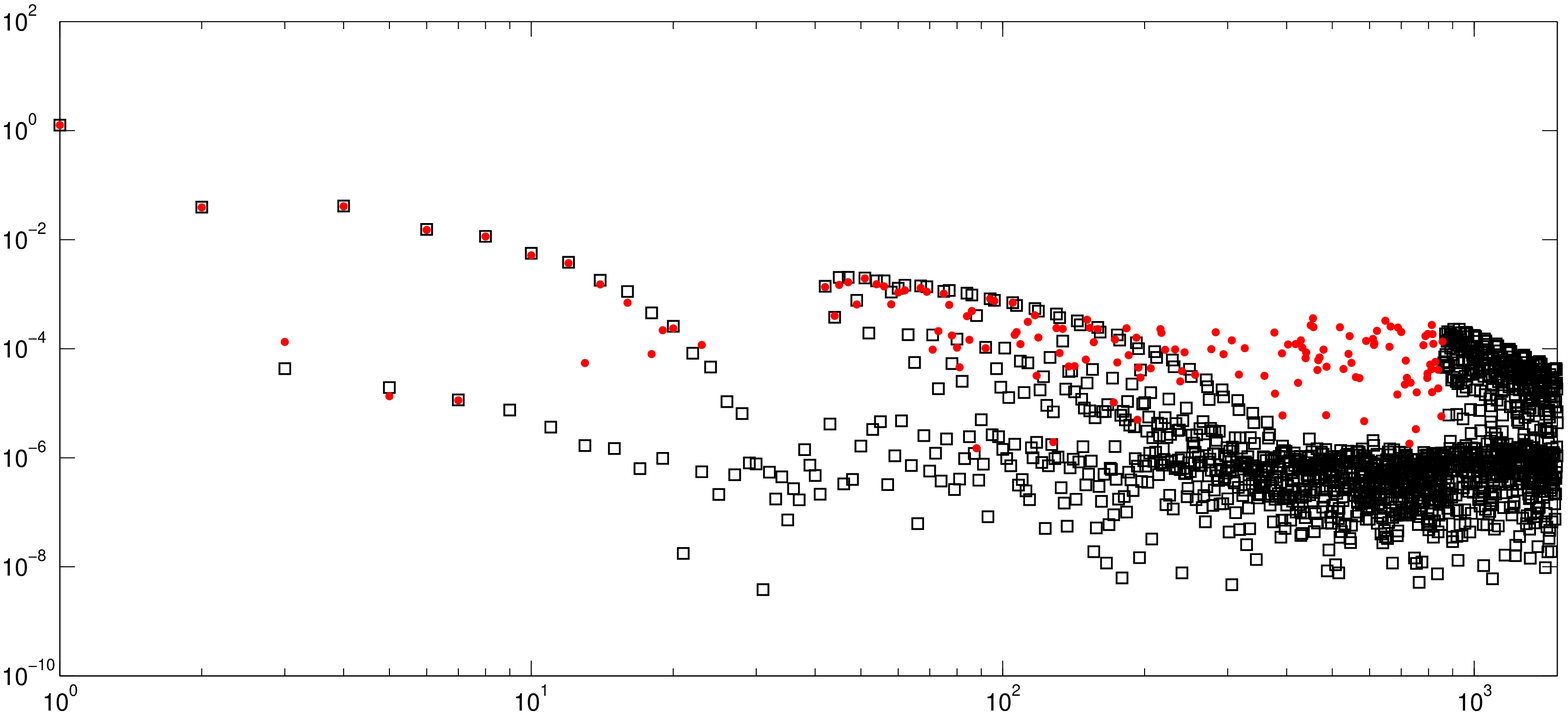}
      \put(-335,60){\begin {sideways} $\vert c_{\bm \alpha_i}\vert$ \end{sideways}}
            \put(-300,25){$\bm{(a)}$}
      \\
      \psfrag{x}{$x$}
      \includegraphics[width=4.8in]{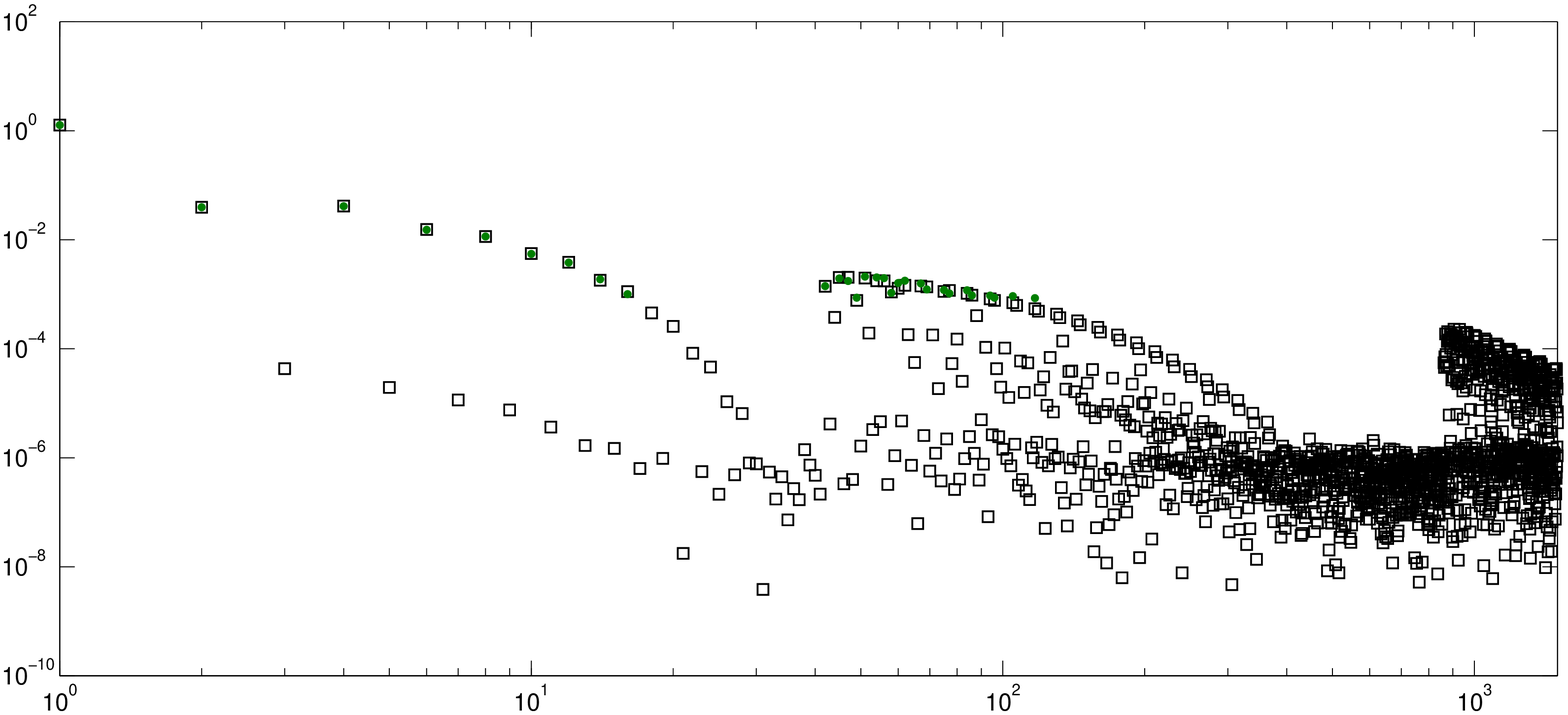}
      \put(-335,60){\begin {sideways} $\vert c_{\bm \alpha_i}\vert$ \end{sideways}}
            \put(-300,25){$\bm{(b)}$}
      \\
      \psfrag{x}{$x$}
      \includegraphics[width=4.8in]{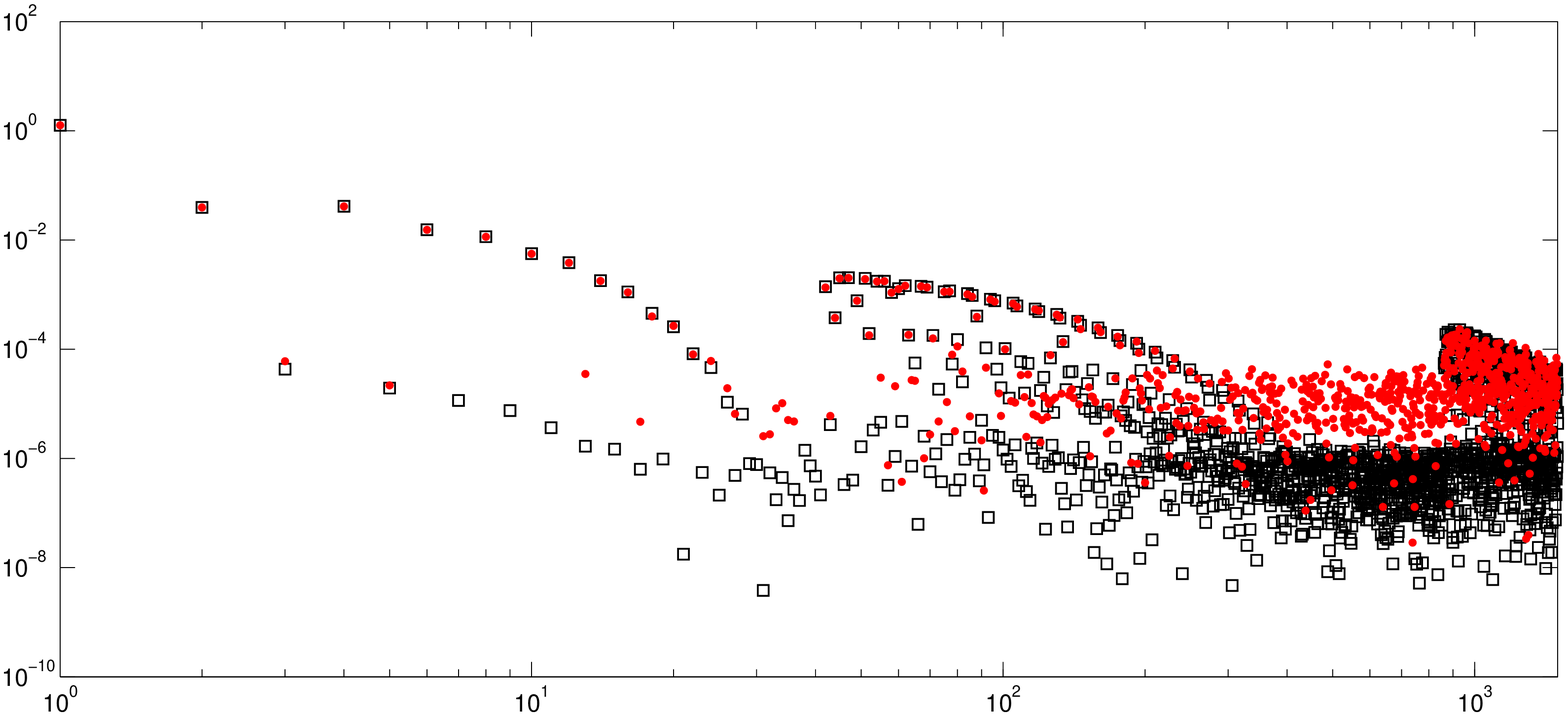}
      \put(-335,60){\begin {sideways} $\vert c_{\bm \alpha_i}\vert$ \end{sideways}}
            \put(-300,25){$\bm{(c)}$}
      \\
      \psfrag{x}{$x$}	
      \includegraphics[width=4.8in]{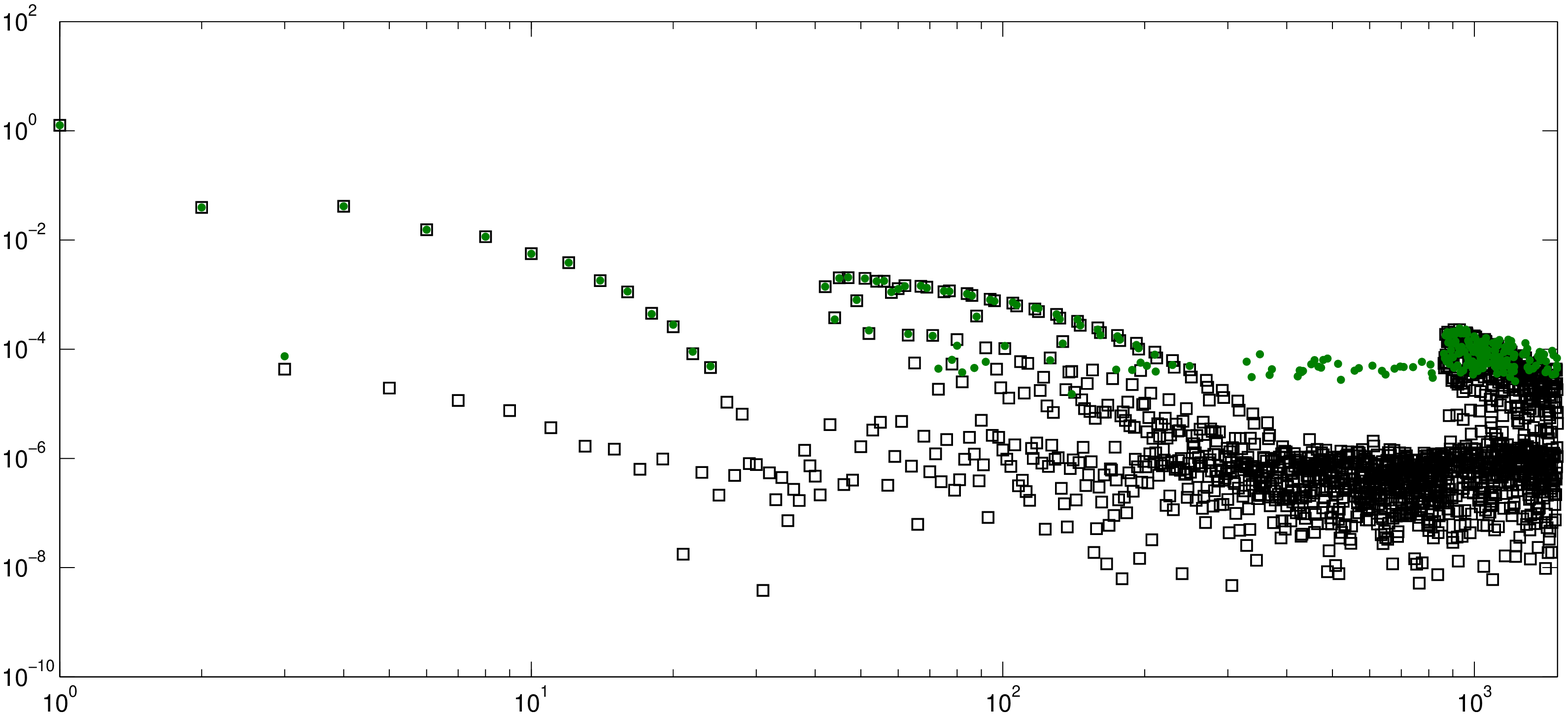}
      \put(-335,60){\begin {sideways} $\vert c_{\bm \alpha_i}\vert$ \end{sideways}}
            \put(-300,25){$\bm{(d)}$}
      \put(-230,-4){\footnotesize Index of PC coefficients ($i$)}
     \end{tabular}
      \caption{Approximation of polynomial chaos (PC) coefficients $\bm c$ of $u(0.5,\bm{y})$ using BPDN and OMP for $d=40$. (a) BPDN with $N=200$ samples, (b) OMP with $N=200$ samples, (c) BPDN with $N=1000$ samples, and (d) OMP with $N=1000$ samples . `Exact' coefficients computed from level 5 stochastic collocation with the Clenshaw-Curtis abscissas ({\scriptsize $\square$}); BPDN and OMP ($\bullet$).}
     \label{fig:pce_coef_40}
\end{figure}
\clearpage

\begin{figure}
    \centering
    \begin{tabular}{cc}

      \psfrag{x}{$x$}	
      \includegraphics[width=2.67in]{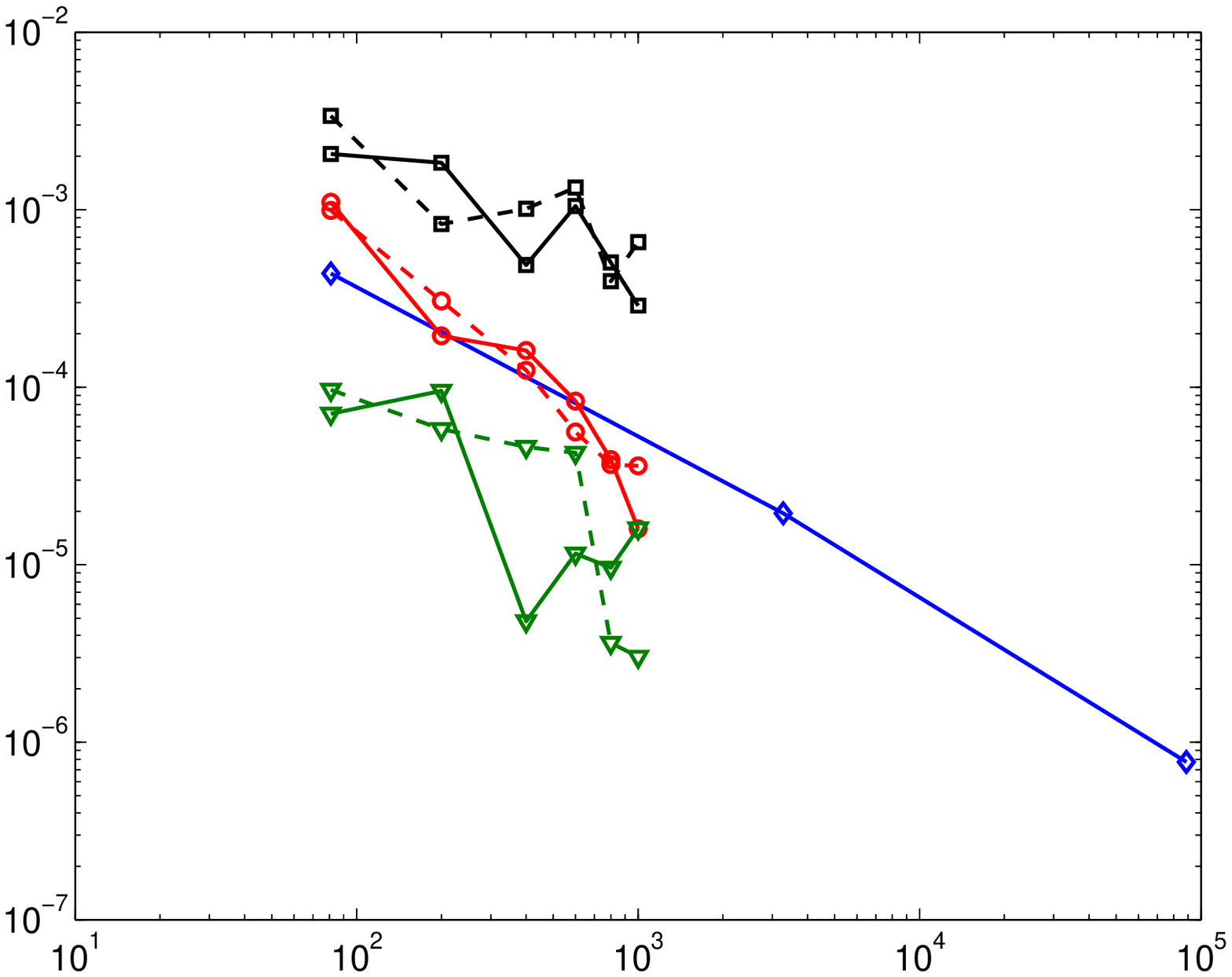}
      \put(-205,35){\begin {sideways} {\footnotesize Rel. error in mean} \end{sideways}}
      \put(-100,-10){$N$}
            \put(-30,130){$\bm{(a)}$}
&
      \psfrag{x}{$x$}	
      \includegraphics[width=2.67in]{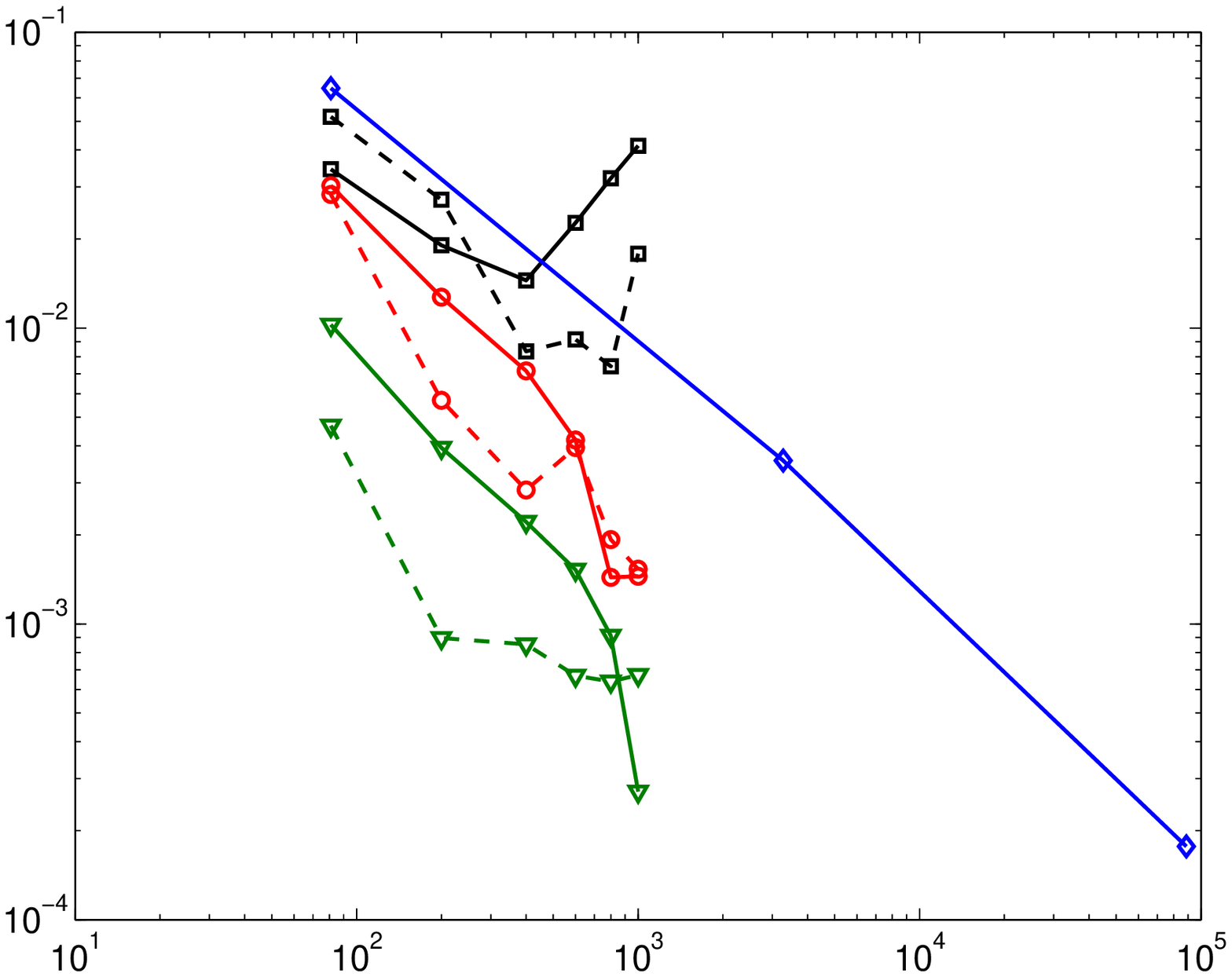}
      \put(-205,35){\begin {sideways} {\footnotesize Rel. error in s.d.} \end{sideways}}
     \put(-100,-10){$N$}
           \put(-30,130){$\bm{(b)}$}
      \\
      \psfrag{x}{$x$}
      \includegraphics[width=2.67in]{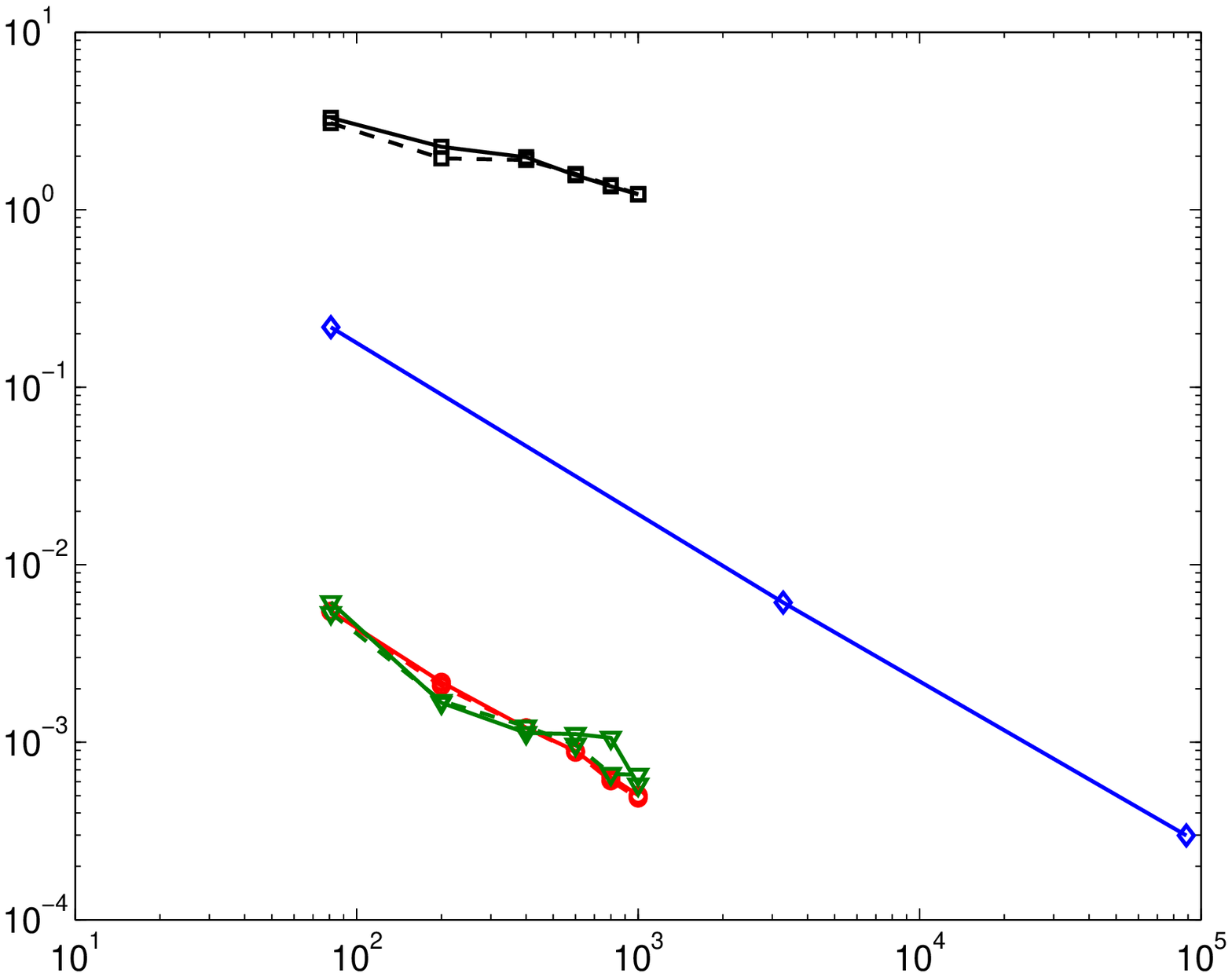}
      \put(-205,40){\begin {sideways} {\footnotesize Rel. rms error} \end{sideways}}
      \put(-100,-10){$N$}
            \put(-30,130){$\bm{(c)}$}
&
      \psfrag{x}{$x$}	
      \includegraphics[width=2.9in]{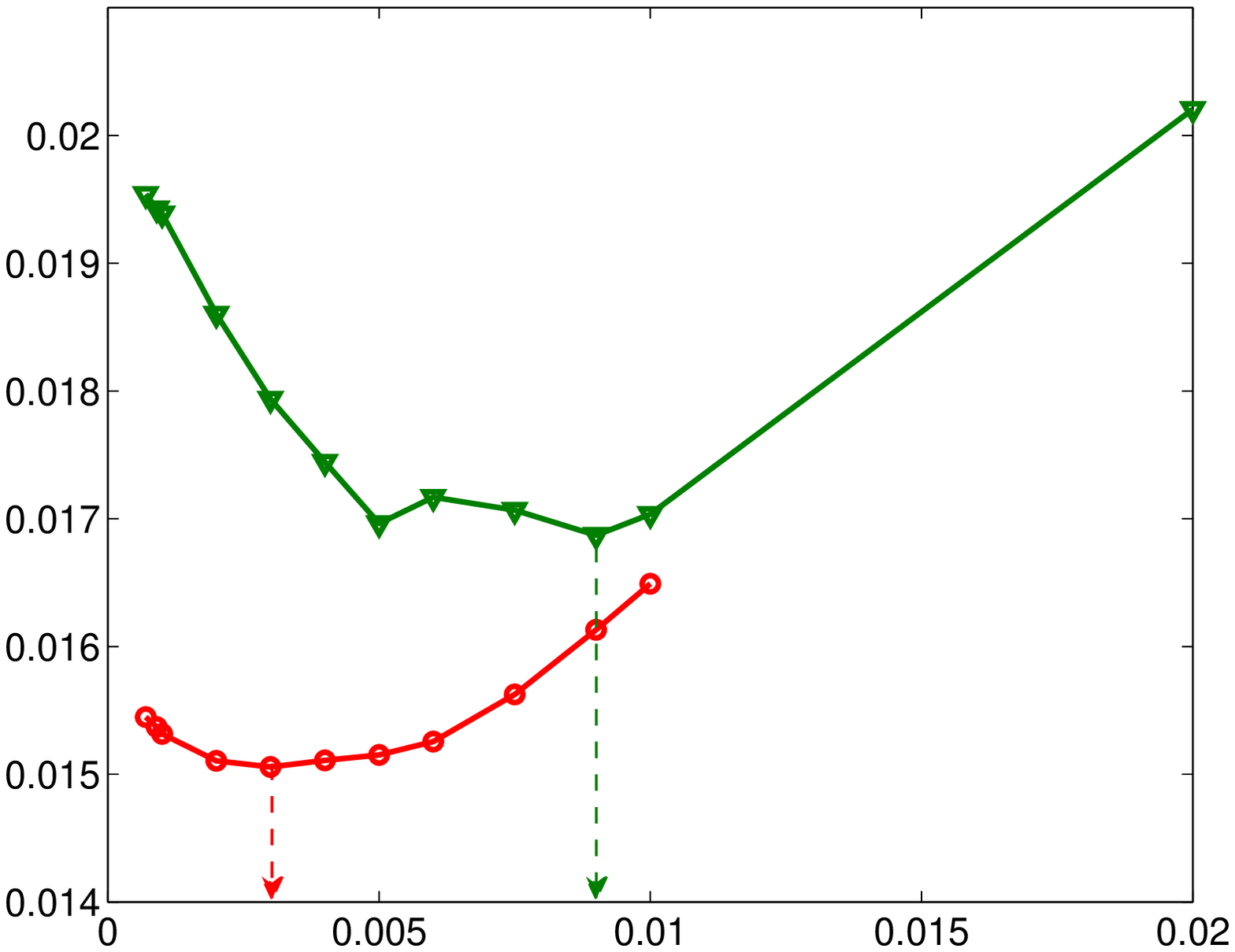}
      \put(-215,75){\begin {sideways} $\delta_{v}$ \end{sideways}}
      \put(-170,140){\footnotesize $N=1000, p=4$}
      \put(-105,0){$\delta_{r}$}
      \put(-115,8){\scriptsize \color{OliveGreen} $\hat{\delta}_{r}$}
      \put(-160,8){\scriptsize \color{red} $\hat{\delta}_{r}$}
            \put(-45,135){$\bm{(d)}$}
     \end{tabular}
          \caption{Comparison of relative error in solution statistics at $x=0.5$ for the Monte Carlo simulation, isotropic sparse grid stochastic collocation with the Clenshaw-Curtis abscissas, and the proposed sparse approximations (BPDN and OMP) for $d=40$. Two sets of independent random samples of $u(0.5,\bm{y})$ are generated first and are used for the Monte Carlo simulation, BPDN, and OMP. The solid and dashed lines correspond to the first and second sets of samples, respectively. (a) Relative error in mean; (b) Relative error in standard deviation; (c) Relative root-mean-squares (rms) error; and (d) Estimation of $\delta$ using cross-validation: $\delta$ is computed from $\delta_{r}$ for which $\delta_{v}$ is minimum. (Monte Carlo simulation ({\scriptsize $\square$}); stochastic collocation ($\diamond$); BPDN ($\circ$); OMP ($\triangledown$)).}
     \label{fig:mean_sd_mse_cross_40}
\end{figure}

For $N=\{81,200\}$ we only estimate the coefficients associated with the $2$nd-order PC expansion, i.e. $p=2$ and $P=861$. For larger sample sizes, we also include the first $639$ basis functions from the $3$th-order chaos, thus leading to $P=1500$. For each combination of $N$ and $p$, we estimate the truncation error $\delta$ using an identical cross-validation procedure described in Section \ref{sec:case_d_40}. Figure \ref{fig:pce_coef_40} illustrates the estimation of PC coefficients of $u(0.5,\bm{y})$ with BPDN and OMP algorithms with $N=200$ and $N=1000$. We again note that the recovered solution from the BPDN algorithm is less sparse as compared to that of the OMP approach, although the over-estimated coefficients (mostly from the second order term) are indeed small. As the samples size $N$ is increased, we are naturally able to recover more dominant coefficients on the expansion. Figure \ref{fig:mean_sd_mse_cross_40} depicts the convergence of the statistics of the solution as functions of the sample size $N$ as well as one instance of the estimation of the truncation error tolerance $\delta$. The implementation details are similar to those described in Section \ref{sec:case_d_14} for the case of $d=14$.

It is worth highlighting that the computational saving of the present sparse approximations (in terms of the number of samples needed to achieve a certain accuracy) compared to the isotropic sparse grid collocation is even larger for the higher-dimensional ($d=40$) problem. This is due to the fact that the number of samples needed to recover the solution using $\ell_1$- and $\ell_0$-minimization is dictated more by the number of dominant terms in the PC expansion compared to the total number of terms $P$, as in Theorem \ref{the:main0}.

\section{Conclusion}
\label{sec:conclusion}
The present study proposes a {\it non-intrusive} and {\it non-adapted} approach based on the {\it compressive sampling} formalism for the approximation of sparse solution of stochastic PDEs. When sufficiently sparse in the polynomial chaos (PC) basis, the compressive sampling enables an accurate recovery of the solution using a number of random solution samples that is significantly smaller than the cardinality of the PC basis. Sparse PC approximations based on $\ell_0$- and $\ell_1$-minimization approaches have been introduced and implemented using the Basis Pursuit Denoising (BPDN) and the Orthogonal Matching Pursuit (OMP) algorithms, respectively. Probabilistic bounds based on the {\it concentration of measure} phenomenon have been derived to verify the convergence and the stability of the present sparse constructions. The performance and efficiency of the proposed techniques are explored through their application to a linear elliptic PDE with high-dimensional random diffusion coefficients where the sparsity of the solution with respect to the PC basis is guaranteed. The proposed formalism to recover the sparse PC expansion of stochastic functions is not restricted to the case of the elliptic PDEs, as its underlying applicability assumptions are universal. Although the discussions of this work have been focused on the particular case of the Legendre PC expansions, the proposed framework can be readily extended to other bases such as the Hermite PC (when the random variables $\bm y$ are standard Gaussian).

\section*{Acknowledgments}
The first author acknowledges the support of the United States Department of Energy under Stanford's Predictive Science Academic Alliance Program (PSAAP) for the preliminary stages of his work. The second author acknowledges the support of the National Science Foundation via NSF grant CMMI-092600 and of the United States Department of Energy under Caltech's Predictive Science Academic Alliance Program (PSAAP).

\bibliographystyle{plain}
\bibliography{AD_bib_v1}

\begin{thebibliography}{10}

\bibitem{Babuska07a}
I.~Babu{\v{s}}ka, F.~Nobile, and R.~Tempone.
\newblock A stochastic collocation method for elliptic partial differential
  equations with random input data.
\newblock {\em SIAM J. Numer. Anal.}, 45(3):1005--1034, 2007.

\bibitem{Babuska04}
I.~Babu{\v{s}}ka, R.~Tempone, and G.~Zouraris.
\newblock Galerkin finite element approximations of stochastic elliptic partial
  differential equations.
\newblock {\em SIAM Journal on Numerical Analysis}, 42(2):800--825, 2004.

\bibitem{Beck09}
A.~Beck and M.~Teboulle.
\newblock A fast iterative shrinkage-threshold algorithm for linear inverse
  problems.
\newblock {\em SIAM J. Imaging Sciences}, 2:183--202, 2009.

\bibitem{Becker09}
S.~Becker, J.~Bobin, and E.~J. Candes.
\newblock {NESTA: A} fast and accurate first-order method for sparse recovery.
\newblock {\em ArXiv e-prints}, 2009.
\newblock Available from http://arxiv.org/abs/0904.3367.

\bibitem{spgl1:2007}
E.~{van den} Berg and M.~P. Friedlander.
\newblock {SPGL1}: A solver for large-scale sparse reconstruction, June 2007.
\newblock Available from http://www.cs.ubc.ca/labs/scl/spgl1.

\bibitem{Bieri09b}
M.~Bieri.
\newblock A sparse composite collocation finite element method for elliptic
  spdes.
\newblock Technical Report Research Report No. 2009-08, Seminar {f\"ur}
  Angewandte Mathematik, SAM, {Z\"urich}, Switzerland, 2009.

\bibitem{Bieri09c}
M.~Bieri, R.~Andreev, and C.~Schwab.
\newblock Sparse tensor discretization of elliptic spdes.
\newblock Technical Report Research Report No. 2009-07, Seminar {f\"ur}
  Angewandte Mathematik, SAM, {Z\"urich}, Switzerland, 2009.

\bibitem{Bieri09a}
M.~Bieri and C.~Schwab.
\newblock Sparse high order {FEM} for elliptic {sPDEs}.
\newblock {\em Computer Methods in Applied Mechanics and Engineering}, 198:1149
  -- 1170, 2009.

\bibitem{Bioucas-Dias07}
J.M. Bioucas-Dias and M.A.T. Figueiredo.
\newblock A new {TwIST}: Two-step iterative shrinking/thresholding algorithms
  for image restoration.
\newblock {\em IEEE Trans. Image Processing}, 16(12):2992--3004, 2007.

\bibitem{Blatman10}
G.~Blatman and B.~Sudret.
\newblock An adaptive algorithm to build up sparse polynomial chaos expansions
  for stochastic finite element analysis.
\newblock {\em Probabilistic Engineering Mechanics}, 25(2):183--197, 2010.

\bibitem{Boufounos07}
P.~Boufounos, M.F. Duarte, and R.G. Baraniuk.
\newblock Sparse signal reconstruction from noisy compressive measurements
  using cross validation.
\newblock In {\em SSP '07: Proceedings of the 2007 IEEE/SP 14th Workshop on
  Statistical Signal Processing}, pages 299--303. IEEE Computer Society, 2007.

\bibitem{Bredies08}
K.~Bredies and D.A. Lorenz.
\newblock Linear convergence of iterative soft- thresholding.
\newblock {\em SIAM J. Sci. Comp.}, 30(2):657--683, 2008.

\bibitem{Bruckstein09}
A.M. Bruckstein, D.L. Donoho, and M.~Elad.
\newblock From sparse solutions of systems of equations to sparse modeling of
  signals and images.
\newblock {\em SIAM Review}, 51(1):34--81, 2009.

\bibitem{Candes06b}
E.J. Candes and J.~Romberg.
\newblock Quantitative robust uncertainty principles and optimally sparse
  decompositions.
\newblock {\em Found. Comput. Math.}, 6(2):227--254, 2006.

\bibitem{Candes07a}
E.J. Candes and J.~Romberg.
\newblock Sparsity and incoherence in compressive sampling.
\newblock {\em Inverse Problems}, 23(3):969--985, 2007.

\bibitem{Candes06a}
E.J. Candes, J.~Romberg, and T.~Tao.
\newblock Robust uncertainty principles: exact signal reconstruction from
  highly incomplete frequency information.
\newblock {\em Information Theory, IEEE Transactions on}, 52(2):489--509, 2006.

\bibitem{Candes06c}
E.J. Candes and T.~Tao.
\newblock Near optimal signal recovery from random projections: Universal
  encoding strategies?
\newblock {\em IEEE Transactions on information theory}, 52(12):5406--5425,
  2006.

\bibitem{Candes08a}
E.J. Candes, M.B. Wakin, and S.~Boyd.
\newblock Enhancing sparsity by reweighted $\ell_1$ minimization.
\newblock {\em Journal of Fourier Analysis and Applications}, 14(5):877--905,
  2008.

\bibitem{Chen98}
S.S. Chen, D.L. Donoho, and M.~Saunders.
\newblock Atomic decomposition by basis pursuit.
\newblock {\em SIAM J. Sci. Comput.}, 20:33--61, 1998.

\bibitem{Chen01a}
S.S. Chen, D.L. Donoho, and M.~Saunders.
\newblock Atomic decomposition by basis pursuit.
\newblock {\em SIAM Rev.}, 43(1):129--159, 2001.

\bibitem{Cohen09a}
A.~Cohen, W.~Dahmen, and R.~DeVore.
\newblock Compressed sensing and best $k-$term approximation.
\newblock {\em J. Amer. Math. Soc.}, 22:211--231, 2009.

\bibitem{Combettes05}
P.L. Combettes and V.R. Wajs.
\newblock Signal recovery by proximal forward-backward splitting.
\newblock {\em Multiscale Modeling and Simulation}, 4(4):1168--1200, 2005.

\bibitem{Dai09}
W.~Dai and O.~Milenkovic.
\newblock Subspace pursuit for compressive sensing signal reconstruction.
\newblock {\em IEEE Trans. Inf. Theor.}, 55(5):2230--2249, 2009.

\bibitem{Daubechies04}
I.~Daubechies, M.~Defrise, and C.~De Mol.
\newblock An iterative thresholding algorithm for linear inverse problems with
  a sparsity constraint.
\newblock {\em Communications on Pure and Applied Mathematics},
  57(11):1413--1457, 2004.

\bibitem{Blumensath09}
T.~Blumensath~M.E. Davies.
\newblock Iterative hard thresholding for compressed sensing.
\newblock {\em Appl. Comp. Harm. Anal}, 27(3):265--274, 2009.

\bibitem{Davis97}
G.~Davis, S.~Mallat, and M.~Avellaneda.
\newblock Adaptive greedy approximation.
\newblock {\em J. Constr. Approx}, 13:57--98, 1997.

\bibitem{Deb01}
M.K. Deb, I.~Babuska, and J.T. Oden.
\newblock Solution of stochastic partial differential equations using
  {G}alerkin finite element techniques.
\newblock {\em Comput. Methods Appl. Mech. Engrg.}, 190:6359--6372, 2001.

\bibitem{Sparselab}
D.~L. Donoho, I.~Drori, V.C. Stodden, and Y.~Tsaig.
\newblock Available from http://www- stat.stanford.edu/{$\sim$}sparselab/.

\bibitem{Donoho06b}
D.L. Donoho.
\newblock Compressed sensing.
\newblock {\em IEEE Transactions on information theory}, 52(4):1289--1306,
  2006.

\bibitem{Donoho06c}
D.L. Donoho, I.~Drori, Y.~Tsaig, and J.L. Starck.
\newblock Sparse solution of underdetermined linear equations by stagewise
  orthogonal matching pursuit.
\newblock Technical report, 2006.

\bibitem{Donoho06a}
D.L. Donoho, M.~Elad, and V.N. Temlyakov.
\newblock Stable recovery of sparse overcomplete representations in the
  presence of noise.
\newblock {\em IEEE Transactions on information theory}, 52(1):6--18, 2006.

\bibitem{Doostan07}
A.~Doostan, R.~Ghanem, and J.~Red-Horse.
\newblock Stochastic model reduction for chaos representations.
\newblock {\em Computer Methods in Applied Mechanics and Engineering},
  196(37-40):3951--3966, 2007.

\bibitem{Doostan09}
A.~Doostan and G.~Iaccarino.
\newblock A least-squares approximation of partial differential equations with
  high-dimensional random inputs.
\newblock {\em Journal of Computational Physics}, 228(12):4332--4345, 2009.

\bibitem{Efron04}
B.~Efron, T.~Hastie, L.~Johnstone, and R.~Tibshirani.
\newblock Least angle regression.
\newblock {\em Annals of Statistics}, 32:407--499, 2004.

\bibitem{Figueiredo07}
M.A.T Figueiredo, R.D. Nowak, and S.J. Wright.
\newblock Gradient projection for sparse reconstruction: Application to
  compressed sensing and other inverse problems.
\newblock {\em Selected Topics in Signal Processing, IEEE Journal of},
  1(4):586--597, 2007.

\bibitem{Ghanem03}
R.~Ghanem and P.~Spanos.
\newblock {\em Stochastic Finite Elements: A Spectral Approach}.
\newblock Dover, 2002.

\bibitem{Hale08}
E.T. Hale, W.~Yin, and Y.~Zhang.
\newblock Fixed-point continuation for $\ell_1$-minimization: Methodology and
  convergence.
\newblock {\em SIAM J. on Optimization}, 19(3):1107--1130, 2008.

\bibitem{Hosder06}
S.~Hosder, R.W. Walters, and R.~Perez.
\newblock A non-intrusive polynomial chaos method for uncertainty propagation
  in {CFD} simulations.
\newblock In {\em $44$th AIAA aerospace sciences meeting and exhibit,
  AIAA-2006-891}, Reno (NV), 2006.

\bibitem{Khajehnejad10}
M.A. Khajehnejad, W.~Xu, A.S. Avestimehr, and B.~Hassibi.
\newblock Improved sparse recovery thresholds with two-step reweighted $\ell_1$
  minimization.
\newblock {\em ArXiv e-prints}, 2010.
\newblock Available from http://arxiv.org/abs/1004.0402.

\bibitem{Kim07a}
S.-J Kim, K.~Koh, M.~Lustig, S.~Boyd, and D.~Gorinevsky.
\newblock An interior-point method for large-scale l1-regularized least
  squares.
\newblock {\em Selected Topics in Signal Processing, IEEE Journal of},
  1(4):606--617, 2007.

\bibitem{Ma09a}
X.~Ma and N.~Zabaras.
\newblock An adaptive hierarchical sparse grid collocation algorithm for the
  solution of stochastic differential equations.
\newblock {\em Journal of Computational Physics}, 228:3084--3113, 2009.

\bibitem{Mathelin03}
L.~Mathelin and M.Y. Hussaini.
\newblock A stochastic collocation algorithm for uncertainty analysis.
\newblock Technical Report NAS 1.26:212153; NASA/CR-2003-212153, NASA Langley
  Research Center, 2003.

\bibitem{Needell09}
D.~Needell.
\newblock Noisy signal recovery via iterative reweighted $l1$-minimization.
\newblock In {\em Proc. Asilomar Conf. on Signals, Systems, and Computers},
  Pacific Grove, CA, Nov. 2009.

\bibitem{Needell08a}
D.~Needell and J.A. Tropp.
\newblock {CoSaMP}: Iterative signal recovery from incomplete and inaccurate
  samples.
\newblock {\em Applied and Computational Harmonic Analysis}, 26(3):301--321,
  2008.

\bibitem{Needell07}
D.~Needell and R.~Vershynin.
\newblock Signal recovery from incomplete and inaccurate measurements via
  {R}egularized {O}rthogonal {M}atching {P}ursuit.
\newblock {\em ArXiv e-prints}, 2007.

\bibitem{Nobile08b}
F.~Nobile, R.~Tempone, and C.G. Webster.
\newblock An anisotropic sparse grid stochastic collocation method for partial
  differential equations with random input data.
\newblock {\em SIAM J. Numer. Anal.}, 46(5):2411--2442, 2008.

\bibitem{Nobile08a}
F.~Nobile, R.~Tempone, and C.G. Webster.
\newblock A sparse grid stochastic collocation method for partial differential
  equations with random input data.
\newblock {\em SIAM J. Numer. Anal.}, 46(5):2309--2345, 2008.

\bibitem{Nouy07}
A.~Nouy.
\newblock {A generalized spectral decomposition technique to solve a class of
  linear stochastic partial differential equations}.
\newblock {\em Computer Methods in Applied Mechanics and Engineering},
  196(37-40):4521--4537, 2007.

\bibitem{Nouy08}
A.~Nouy.
\newblock {Generalized spectral decomposition method for solving stochastic
  finite element equations: Invariant subspace problem and dedicated
  algorithms}.
\newblock {\em Computer Methods in Applied Mechanics and Engineering},
  197:4718--4736, 2008.

\bibitem{Osborne00}
M.R. Osborne, B.~Presnell, and B.~Turlach.
\newblock A new approach to variable selection in least squares problems.
\newblock {\em IMA J. Numer. Anal.}, 20:389--403, 2000.

\bibitem{Pati93}
Y.~C. Pati, R.~Rezaiifar, and P.~S. Krishnaprasad.
\newblock {O}rthogonal {M}atching {P}ursuit: {R}ecursive function approximation
  with applications to wavelet decomposition.
\newblock In {\em Proceedings of the 27th Annual Asilomar Conference on
  Signals, Systems, and Computers}, pages 40--44, 1993.

\bibitem{Schwab06a}
C.~Schwab and R.A. Todor.
\newblock Karhunen-{L}o\`{e}ve approximation of random fields by generalized
  fast multipole methods.
\newblock {\em J. Comput. Phys.}, 217(1):100--122, 2006.

\bibitem{Tatang95}
M.A. Tatang.
\newblock {\em Direct incorporation of uncertainty in chemical and
  environmental engineering systems}.
\newblock PhD thesis, Ph.D. Thesis, Massachusetts Institute of Technology,
  1995.

\bibitem{Tibshirani94}
R.~Tibshirani.
\newblock Regression shrinkage and selection via the {L}asso.
\newblock {\em Journal of the Royal Statistical Society, Series B},
  58(1):267--288, 1996.

\bibitem{Todor07a}
R.~A. Todor and C.~Schwab.
\newblock Convergence rates for sparse chaos approximations of elliptic
  problems with stochastic coefficients.
\newblock {\em IMA Journal of Numerical Analysis}, 27(2):232--261, 2007.

\bibitem{Tropp10a}
J.A. Tropp and S.J. Wright.
\newblock Computational methods for sparse solution of linear inverse problems.
\newblock {\em Proceedings of the IEEE}, 2010.
\newblock in press.

\bibitem{Berg08}
E.~van~den Berg and M.~P. Friedlander.
\newblock Probing the {P}areto frontier for basis pursuit solutions.
\newblock {\em SIAM Journal on Scientific Computing}, 31(2):890--912, 2008.

\bibitem{Wan05}
X.~Wan and G.~Karniadakis.
\newblock An adaptive multi-element generalized polynomial chaos method for
  stochastic differential equations.
\newblock {\em J. Comp. Phys.}, 209:617--642, 2005.

\bibitem{Ward09}
R.~Ward.
\newblock Compressed sensing with cross validation.
\newblock {\em IEEE Trans. Inf. Theor.}, 55(12):5773--5782, 2009.

\bibitem{Xiu05a}
D.~Xiu and J.S. Hesthaven.
\newblock High-order collocation methods for differential equations with random
  inputs.
\newblock {\em SIAM J. Sci. Comput.}, 27(3):1118--1139, 2005.

\bibitem{Xiu02}
D.~Xiu and G.M. Karniadakis.
\newblock The {W}iener-{A}skey polynomial chaos for stochastic differential
  equations.
\newblock {\em SIAM Joural on Scientific Computin}, 24(2):619--644, 2002.

\bibitem{Xu09a}
W.~Xu, M.A. Khajehnejad, A.S. Avestimehr, and B.~Hassibi.
\newblock Breaking through the thresholds: an analysis for iterative reweighted
  $\ell_1$ minimization via the {G}rassmann {A}ngle {F}ramework.
\newblock {\em ArXiv e-prints}, 2009.
\newblock Available from http://arxiv.org/abs/0904.0994.

\bibitem{Yang09}
J.~Yang and Y.~Zhang.
\newblock Alternating direction algorithms for $\ell_1$-problems in compressive
  sensing.
\newblock {\em ArXiv e-prints}, 2009.
\newblock Available from http://arxiv.org/abs/0912.1185.

\end{thebibliography}

\end{document}